\documentclass[a4paper,reqno,english]{smfart}

\usepackage{amssymb,amsthm,xspace,smfthm}
\usepackage{latexsym}
\usepackage{natbib} 
\usepackage{lineno}
\usepackage{amsmath}
\usepackage{mathrsfs}
\usepackage{epsfig}
\usepackage[latin1]{inputenc}
\usepackage[normalem]{ulem}
\usepackage{verbatim}
\usepackage{graphicx}
\usepackage{graphics}
\usepackage{color}
\usepackage{enumerate}
\usepackage{bbm}
\usepackage{dsfont}
\usepackage{bm}
\usepackage{graphics}
\usepackage{epsfig}
\usepackage{fleqn}
\usepackage{appendix}

\newcommand{\XX}{\Xi}
\newcommand{\p}{p}
\newcommand{\Vol}{\mbox{Vol}}

\newcommand{\U}{\mathbbm{U}}
\newcommand{\Ss}{\mathcal{S}}
\newcommand{\cal}{\mathcal}

\newcommand{\E}{\mbox{E}}
\newcommand{\F}{\mathbbm{F}}

\newcommand{\Var}{\mbox{Var}}

\newcommand{\1}{\mathbbm{1}}
\newcommand{\s}{\xi}

\newcommand{\Ll}{\mathbbm{L}}

\newcommand{\R}{I\!\!R}

\newcommand{\Pp}{I\!\!P}

\theoremstyle{remark}

\newtheorem{lemma}{Lemma}
\newtheorem{deff}{Definition}

\newtheorem{proposition}{Proposition}[section]

\newtheorem{theorem}{Theorem}[section]

\theoremstyle{definition}
\newtheorem{remark}{Remark}

\newcommand*{\EB}[1]{%
  \makebox[0pt][l]{%
    \kern.1em
    \raisebox{1.55ex}{\relsize{-3}{$\circ$}}
  }%
  #1%
}

\begin{document}
\title[Accelerated Monte Carlo under monotonicity]{Accelerated Monte Carlo estimation of exceedance probabilities under monotonicity constraints}

\author[Nicolas Bousquet]{Nicolas Bousquet}

\maketitle

\begin{center}
\footnotesize 
\vspace{-1cm}
{\it EDF Research \& Development} \\
{\it Dpt. of Industrial Risk Management} \\
{\it 6 quai Watier, 78401 Chatou, France} \\
\texttt{nicolas.bousquet@edf.fr}
\vspace{1cm}
\normalsize
\end{center}

\begin{abstract}
The problem of estimating the probability $\p=P(g({\bf X})\leq 0)$ is considered when ${\bf X}$ represents a multivariate stochastic input of a monotonic function $g$. First, a heuristic method to bound $\p$, originally proposed by \citet{ROC09}, is  formally described, involving a specialized design of numerical experiments. Then a statistical estimation of $\p$ is considered based on a sequential stochastic exploration of the input space. A maximum likelihood estimator of $\p$ based on successive dependent Bernoulli data is defined and its theoretical convergence properties are studied. Under intuitive or mild conditions, the estimation is faster and more robust than the traditional Monte Carlo approach, therefore adapted to time-consuming computer codes $g$. The main result of the paper is related to the variance of the estimator. It appears as a new baseline measure of efficiency under monotonicity constraints, which could play a similar role to the usual Monte Carlo estimator variance in unconstrained frameworks. Furthermore the bias of the estimator is shown to be corrigible via bootstrap heuristics. The behavior of the method is illustrated by numerical tests led on a class of toy examples and a more realistic hydraulic case-study. \\

 On considère l'estimation de la probabilité $\p=P(g({\bf X})\leq 0)$ où $\bf X$ est un vecteur aléatoire  et $g$ une fonction monotone. Premièrement, on rappelle et formalise une méthode, proposée par de Rocquigny (2009), permettant d'encadrer $\p$ par des bornes déterministes en fonction d'un plan d'expérience séquentiel. Le second et principal apport de l'article est la définition et l'étude  d'un estimateur statistique de $\p$ tirant parti des bornes. Construit à partir de tirages uniformes successifs, cet estimateur présente sous de faibles conditions théoriques une variance asymptotique plus faible et une meilleure robustesse que l'estimateur classique de Monte Carlo, ce qui rend la méthode adaptée à l'emploi de codes informatiques $g$ lourds en temps de calcul.  Des expérimentations numériques sont menées sur des exemples-jouets et un cas d'étude hydraulique plus réaliste. Une heuristique de boostrap, reposant sur un réplicat de l'hypersurface $\{{\bf x}, \ g({\bf x})= 0\}$ par des réseaux de neurones, est proposée et testée avec succès pour ôter le biais non-asymptotique de l'estimateur. 
\end{abstract}

%
%
%
%

\section{Introduction}\label{introduction}

In many technical areas, the exceedance of some unidimensional variable $Z$ over a certain critical value $z^*$ 
may define  an event of probability $\p$ which has to be carefully estimated. Assumed to be stricly positive, $\p$  can be defined by 
\begin{eqnarray}
\nonumber \p & = & P(g({\bf X})\leq 0) \ = \  \int_{\U} \1_{\{g({\bf{x}})\leq 0\}}f({\bf{x}}) \ d{\bf{x}} \label{uncertainty.model}
\end{eqnarray}
with $\bf X$ a random vector of uncertain input parameters with probability density function (pdf) $f$, taking its values in a $d-$dimensional space $\U$, and $g({\bf X})=z^*-Z$  a deterministic mapping from $\U$ to $\R$. This framework is often encountered in  structural reliability studies \citep{MAD96}, when $g$ is a computer code reproducing a physical phenomenon. 
 A Monte Carlo (MC) method is the usual way to estimate $\p$, by $\hat{\p}_n=n^{-1}\sum_{k=1}^n \1_{\{g({\bf x}_k)\leq 0\}}$ with large $n$ ${\bf x}_k$ independently sampled from $f$. Avoiding regularity hypotheses on $g$, this unbiased estimator presents good convergence properties  
and an estimation error independent on $d$. Unfortunately, this strategy often appears inappropriate in practice when $\p$ reaches low values, since $g$ can be time-consuming and the computational budget may be limited:  a good estimation of a probability $\p\sim 10^{-q}$ typically requires at least $10^{q+2}$ calls to $g$ \citep{LEM06}. Furthermore, $\hat{\p}_n$ has the theoretical defect not to be robust, in the sense given by \citet{GLY09}: its relative error, namely its coefficient of variation, does not tend to a finite limit when $\p\to 0^+$, given any finite number $n$ of trials. 

Many non-intrusive strategies have been proposed to accelerate the MC approach. Traditional methods from the engineer community in structural reliability (FORM/SORM) treat the estimation of $\p$ as an optimization 
 problem. The computational work is usually fast but the estimators suffer from weakly or non-controllable errors. Statistical approaches are judged in terms of reduction rate with respect to the MC estimator variance $\Var[\hat{p}_n]=\p(1-\p)/n$. Methods like quasi-MC, sequential MC or importance sampling \citep{KRO07} are based on selecting a {\it design of experiments} (DOE), namely a set of points in $\U$ on which $g$ is tested, such that $\U$ be explored in areas close to the limit state surface $\Ss=\{{\bf x}\in\U \ ; \ g({\bf x})=0\}$.  Most advanced methods often get rid of the time-consuming difficulties by emulating the behavior of $g$, for instance using kriging techniques \citep{CAN08} which presuppose smoothness conditions on $g$. 

Minimizing the strength of regularity hypotheses placed on $g$ underlies the development of specialized acceleration methods. For instance, computer codes can suffer from edge effects which restrict smoothness conditions \citep{MUN10}. On the other hand, the reality of the phenomenon can imply various form constraints on $Z$. Especially, the assumption that  $g$ is monotonic with respect to ${\bf X}$ is a standard problem in regression analysis \citep{DUROT08}. In the area of numerical experiments, monotonicity properties of computer codes have been considered theoretically and practically, e.g. proving the MC acceleration of Latin Hypercube Sampling   for the estimation of expectancies \citep{MAC79}, carrying out screening methods for sensitivity analyses \citep{LIN93}, constraining response surfaces  \citep{KLE09,KLE11}, predicting the behavior of network queuing systems \citep{RAN08}, computing flood probabilities \citep{ROC09} or estimating the safety of a nuclear reactor pressure vessel \citep{MUN10}.

Specific engineering works in structural reliability have highlighted the possibility of bounding and estimating $\p$ significantly faster than using a MC approach.  Under the name of {\it monotonic reliability methods} (MRM), \citet{ROC09} proposed a class of sequential  algorithms contouring the limit state surface and enclosing $\p$ between deterministic bounds which dynamically narrow. A similar idea was explored by \citet{RAJ10}. However, although a parallelization of such algorithms was already implemented \citep{LIM10}, these methods were only empirically studied and some of the proposed estimators of $\p$ remained crude. \\

The present article therefore aims to provide a first theoretical approach of the accelerated MC estimation of $\p$ when $g$ is assumed to be monotonic and possibly discontinuous, although some smoothness constraints are assumed on the failure surface $\Ss$. More precisely, this article is structured as follows. 

Section 2 is dedicated to a general description and a mathematical formalization of MRM. The main contribution is presented in Section 3: a statistical estimator of $\p$ is proposed, based on uniformly sampled DOEs in nested spaces. Defined as the maximum likelihood estimator of dependent Bernoulli data, its asymptotic properties are theoretically studied. The estimator is shown to be robust, and its variance gains a significant reduction with respect to the usual MC case. It may also be viewed as a baseline (or {\it target}) variance for monotonic structural reliability frameworks. The non-asymptotic bias of the estimator is examined in Section 4, through numerical experiments involving a class of toy examples. Based on a neural network emulation of $\Ss$, bootstrap heuristics are proposed and successfully tested to remove this bias. Finally, a more realistic hydraulic case-study illustrates the benefits of the complete method. 

Along the paper some connections are done with other areas of computational mathematics, especially about implementation issues, and a discussion section ends this article by focusing on the research avenues that must be explored in the area of stochastic sequential DOEs to improve the results presented here.

\section{Material}
\label{material}

\Subsection{Working assumptions, definitions and basic properties}

Let $g:{\bf X}\mapsto g({\bf X})$ be a deterministic function defined as a
real-valued scalar mapping of ${\bf{X}} = (X_1,\ldots,X_d)$ on
its
definition domain $\U\subset \R^d$.  {\it Deterministic} means that the function $g({\bf x})$ produces always the same output if it is given the same input $\bf x$. {Global monotonicity} is defined as follows:  $\forall i$, $\exists s_i \in\{-1,+1\}$, $\forall \epsilon >0$,  $\forall {\bf x}=(x_1,\ldots,x_d)\in \U$, such that
\begin{eqnarray*}
 g\left(x_1,\ldots,x_{i-1},x_i+s_i\epsilon,x_{i+1},\ldots,x_d\right) & \leq & \hspace{0.25cm} g\left(x_1,\ldots,x_{i-1},x_i,x_{i+1},\ldots,x_d\right)
\end{eqnarray*}
where $s_i$ represents the sign of monotonic dependence: $s_i = 1$ ({\it resp.} $s_i = -1$) when $g$ is decreasing ({\it resp.} increasing) along with the $i-$th component $x_i$. The following assumption is made without loss of generality since any
decreasing $i-$th component can be changed from $x_i$ to $-x_i$: 

\paragraph*{\bf Assumption 1}
The function $g$ is  globally increasing over $\U$. \\

To be general, $\U=[0,1]^d$ and ${\bf X}$ is a random vector defined on the probability space  $(\U,{\cal{B}}(\U),P)$. Next assumption is made following this same concern of generality. 

\paragraph*{\bf Assumption 2}
All inputs $x_1,\ldots,x_d$ are independently uniform on $\U=[0,1]^d$. \\

In real cases, $x_1,\ldots,x_d$ can be defined as {\it transformed inputs}, as usual in structural safety problems \citep{MAD96}. In such cases one can write ${\bf x}=T({\bf y})$ where ${\bf y}=(y_1,\ldots,y_d)$ is a vector of {\it physical} inputs and $T$ is the multivariate distributional transform \citep{RUS09}.  Therefore $g=\tilde{g}\circ T^{-1}$ where $\tilde{g}$ is a mononotic function and $T$ has to preserve this monotonicity. When the $y_i$ are independent, $T$ is reduced to the vector of marginal cdfs $(F_1,\ldots,F_d)$ and is naturally increasing, so the assumption is not restrictive. Else, technical requirements on $T$ are needed, which depend on the way this joint distribution is defined \citep{RUS09}. 
 See for instance \citet{CHEN09} for such requirements on Gaussian copulas. Another general result is given in the Appendix (Supplementary Material). 

\paragraph*{\bf Assumption 3} Both subspaces $\U^-=\{{\bf x}\in\U, \ g({\bf x})\leq 0\}$ and $\U^+=\{{\bf x}\in\U, \ g({\bf x})>0\}$ are not empty (so that $\p$ exists in $]0,1[$).

\begin{deff}\label{domination}
A set of points of $\U$ is said to be {\it safety-dominated} ({\it resp.} {\it failure-dominated}) if $g$ is guaranteed to be positive ({\it resp.} negative) in any point of this set. 
\end{deff}

\noindent Denote by $\succeq$ the partial order
between elements of  $\U$ defined by ${\bf x}\succeq {\bf y} \Leftrightarrow x_k\geq y_k$ $\forall k=1,\ldots,d$. 
 Then
assume that some point value $g({\bf \tilde{x}})$ is known,
and consider the sets 
$\U^+_{{\bf \tilde{x}}}  =  \left\{{\bf x}\in\U \ | \ {\bf x}\succeq {\bf \tilde{x}} \right\}$ and  $\U^-_{{\bf \tilde{x}}}  =  \left\{{\bf x}\in\U \ | \ {\bf x}\preceq {\bf \tilde{x}} \right\}$.
The increasing monotonicity
implies that if $g({\bf \tilde{x}})>0$ ({\it resp.} $g({\bf \tilde{x}})<0$),
then $\U^+_{{\bf \tilde{x}}}$ is { safety-dominated} ({\it resp.} $\U^-_{{\bf \tilde{x}}}$ is { failure-dominated}). This proves next lemma. 

\begin{lemma}\label{lemma.1}
Both inequalities are true with probability 1:
\begin{eqnarray*}
\p & \leq &  1 - P({\bf X}\in\U^+_{{\bf \tilde{x}}})  \ \ \ \text{if $g({\bf \tilde{x}})>0$,} \\
\p & \geq &  P({\bf X}\in\U^-_{{\bf \tilde{x}}}) \ \ \ \text{else.} 
\end{eqnarray*}
\end{lemma}

\noindent More generally, assume that $n$ input vectors $({\bf x}_j)_{j = 1,\ldots, n}$ can be sorted into safe and failure sub-samples following the corresponding values of $\{g( {\bf x}_j )\}_{j = 1,\ldots, n}$. They are respectively defined by
\begin{eqnarray*}
\XX^+_n  & =  & \left\{{\bf x}\in ({\bf x}_j)_{j = 1,\ldots, n} \ | \ g( {\bf x}_j )> 0 \right\}
\end{eqnarray*}
 and 
\begin{eqnarray*}
\XX^-_n  =  \left\{{\bf x}\in ({\bf x}_j)_{j = 1,\ldots, n} \ | \ g( {\bf x}_j )\leq 0 \right\}.
\end{eqnarray*}
Then one may define the sets 
\begin{eqnarray*}
\U^+_n & = & \left\{{\bf x}\in \U \ | \ \exists {\bf x}_j\in \XX^+_n, \ {\bf x}\succeq {\bf x}_j \right\}, \\
 \U^-_n & = & \left\{{\bf x}\in \U \ | \ \exists {\bf x}_j\in \XX^-_n, \ {\bf x}\preceq {\bf x}_j \right\} 
\end{eqnarray*}
(see Figure \ref{graphic.example.2} for an illustration). Finally, denoting $p^-_n=P({\bf X}\in \U^-_n)$ and $p^+_n=1-P({\bf X}\in \U^+_n)$ to alleviate the notations, one has in all the sequel and for all $n\geq 0$,
\begin{eqnarray}
p^-_n & \leq \ \p \ \leq & p^+_n. \label{gen.bounds} 
\end{eqnarray}

Hereafter, $\U^+_n$ and $\U^-_n$  will be referred to as {\it
dominated subspaces}, where the sign of $g({\bf x})$ is known. Note that
the complementary  {\it non-dominated} subspace
$\U_n  = \U / \left(\U^+_n\cup\U^-_n\right)$ 
is the only partition of $\U$ where further calls of $g$ are required to improve the bounds. Finally, a topological assumption on $\Ss$ is  needed to complete  the formal description of the situations studied by \citet{ROC09} and \citet{LIM10}. 

\paragraph*{\bf Assumption 4} The limit state surface $\Ss=\{{\bf x}\in\U \ ; \ g({\bf x})=0\}$ is regular enough and separates $\U$ in two disjoint domains $\U^-$ and $\U^+$ ({\it simply connected}). \\

The second part of this assumption implies that, in terms of {classification}, the two classes of points $\U^-_n$ and $\U^+_n$ are perfectly separable when $n\to\infty$. This property will be used later in the paper to carry out bootstrap heuristics.    
By {\it regular enough}, $\Ss$ is assumed not to be the surface of multidimensional stairs, so that it cannot be exhaustively described by a $n-$DOE with $n<\infty$. This  mild assumption is ensured, for instance, if  $g$ is continuously differentiable on a non-empty measurable subset of $\Ss$.  More formally, it is assumed that $\forall n<\infty$, 
\begin{eqnarray}
\sup\limits_{{\bf x_n}\in\bar{\Ss}} \int_{\U_{n-1}\cap\U^-} \1_{\{{\bf x}\preceq {\bf x_n}\}} \ d{\bf x} & < & \p - \p^-_{n-1}, \label{cond001} \\
\sup\limits_{{\bf x_n}\in\bar{\Ss}} \int_{\U_{n-1}\cap\U^+} \1_{\{1-{\bf x}\preceq 1-{\bf x_n}\}} \ d{\bf x} & < & \p^+_{n-1}-\p. \label{cond002}
\end{eqnarray}
This will imply that $\p^-_n<\p<\p^+_n$  and the finiteness of the strictly positive quantity $\tilde{\omega}_{n+1}(\p)=[(p^{+}_{n}-\p)(\p-p^-_{n})]^{-1}$ encountered further in the paper.

\begin{remark}
In multi-objective optimization, a dominated space can be interpreted as a subset of a performance space 
delimited by a Pareto frontier \citep{FIG05}. In this framework, $g$ is thought as a monotonic rule of decision 
 depending of $d$ variables, for which the set of $n$ best possible configurations (the frontier) is searched. 
\end{remark}  

\begin{remark}
The proportions $(p^-_n,1-p^+_n)$ are the volumes of two unions of hyperrectangles sharing the same orthogonal basis. Computing such volumes is known in computational geometry as Klee's measure problem,  for which recursive sweepline algorithms \citep{LEW81} can provide exact solutions.  Details about their implementation are given in Appendix (Supplementary Material). When $d$ exceeds 4 or 5, these exact methods appear however too costly, and trivial MC methods must be preferred in practice to compute these quantities. 
\end{remark}

\Subsection{MRM implementation: a one-step ahead strategy}

Starting from $\U^+_0=\{1^d\}$, $U^-_0=\{0^d\}$ and $\U_0=\U=[0,1]^d$, the iterative scheme shared by all MRM variants at step $n\geq 1$ is based on: \\

\hbox{\raisebox{0.4em}{\vrule depth 0pt height 0.4pt width 10.5cm} }
\begin{enumerate}
\item selecting a DOE $\{{\bf x}^{(1)}_{n},\ldots,{\bf 
x}^{(m_n)}_n\}\in\U_{n-1}$;
\item computing the {\it signatures} 
$\s^{(j)}_{\bf x_n}  =  \1_{\left\{g\left(\bf x^{(j)}_n\right)<0\right\}}$; 
\item updating the subspaces 
\begin{eqnarray*}
\hspace{2cm} {\U}^-_n  & = &  \U^-_{n-1} \cup \left\{{\bf x}\in\U \ | \ \exists \ {\bf x}^{(j)}_n, \  \ \s^{(j)}_{\bf x_n}=1, \ \ {\bf x} \preceq  {\bf x}^{(j)}_n \right\}, \\
\hspace{2cm} {\U}^+_n  & =  & \U^+_{n-1} \cup \left\{ {\bf x}\in\U \ | \ \exists \
{\bf x}^{(j)}_n, \ \ \s^{(j)}_{\bf x_n}=0, \ \ {\bf x} \succeq
{\bf x}^{(j)}_n\right\}, \\
\hspace{2cm}  \U_{n}  & = &  \U/(\U^{-}_{n}\cup\U^+_{n})
\end{eqnarray*}
\item updating the bounds
$\{p^-_{n},p^+_{n}\}=\{\Vol({\U}^-_n),1-\Vol({\U}^+_n)\}$.
\end{enumerate}
\hbox{\raisebox{0.4em}{\vrule depth 0pt height 0.4pt width 10.5cm} }
\vspace{0.25cm}

Since $\U^-_n\subset\U^-_{n+1}$ $\forall n\geq 0$, then $P({\bf X}\in \U^-_n)\leq P({\bf X}\in \U^-_{n+1})$ and the sequence $(p^-_n)$ is nondecreasing. Symmetrically, the sequence $(p^+_n)$ is nonincreasing. Since bounded in $[0,\p]$ and $[\p,1]$, both sequences are converging.

At each step, the DOE must be chosen accounting for the
increasing monotonicity of $g$. Denoting
$\bf x^{(1)}_n$ and $\bf x^{(2)}_n$ two elements of the DOE and assuming
to know $\s^{(1)}_{\bf x_n}$, it is
unnecessary to compute $\s^{(2)}_{\bf x_n}$ in two cases:
\begin{eqnarray*}
 \text{if $\s^{(1)}_{\bf x_n}=1$ and $\bf x^{(1)}_n \succeq {\bf x}^{(2)}_n
$} & \Rightarrow &  \text{$\bf x^{(2)}_n\in\U^-_{\bf x^{(1)}_n}$ and
$\s^{(2)}_{\bf x_n}=1$,} \label{dom1}\\
 \text{if $\s^{(1)}_{\bf x_n}=0$ and $\bf x^{(1)}_n
\preceq {\bf x}^{(2)}_n $} &  \Rightarrow &
\text{$\bf x^{(2)}_n\in\U^+_{\bf x^{(1)}_n}$ and $\s^{(2)}_{\bf x_n}=0$.}\label{dom2}
\end{eqnarray*}
Thus the order of trials 
 should be carefully monitored, in relation with the partial order between the elements of the DOE. Reducing the DOE to a single element, i.e. $m_n=1$ for all steps, minimizes the number of unnecessary trials. This one-step ahead strategy is favored in the present paper.

\Subsection{Stochastic MRM}

\paragraph*{Initialization} First iterations should be monitored to reduce significantly the width of $[p^-_n,p^+_n]$, such that further iterations mainly focus on refinements. A deterministic strategy seems the most appropriate to start from $[0,1]$ until providing non-trivial bounds. A dichotomic diagonal MRM, illustrated on Figure \ref{types.of.surfaces} in a two-dimensional case, was used in the examples considered further. It explores the non-dominated space in an intuitive way and stops at step $k_0\geq 1$ such that
\begin{eqnarray*}
k_0 & \geq & 1+\frac{\log(1/\p)}{d\log 2}.
\end{eqnarray*}
Consequently, an expected crude prior value of $\p$ can help to estimate the minimal number $k_0$ of  trials. To alleviate the paper, the notation $(\U^+_0,\U^-_0,p^+_0,p^-_0)$ now describes the situation after $N-1$ introductive deterministic  steps with $N\geq k_0+1$, such that $0<p^-_0$ and $p^+_0<1$. 

\paragraph*{Switching to stochastic DOEs} Pursuing a deterministic strategy can be too costly to be efficient, the upper bound $p^+_n$ offering possibly a very conservative assessment of $\p$  \citep{ROC09}. Intuitively, such a strategy should be optimized by selecting the next element of the DOE as the maximizer of a criterion which predicts a measure of dominated volume. Apart from the difficulty of predicting, choosing the criterion remains arbitrary. Switching to a stochastic strategy, which allows for a sequential statistical estimation of $\p$ in addition of providing bounds, seems a promising alternative approach. In this framework, 
\begin{eqnarray*}
{\bf x_n} & \sim & f_{n-1}
\end{eqnarray*} 
at each step $n\geq 1$, with $f_{n-1}$ a pdf defined on $\U_{n-1}$. Then the probability space $(\U,{\cal{B}}(\U),P)$ becomes endowed with the filtration $\F=({\cal{F}}_n)$ where ${\cal{F}}_n$ is the $\sigma-$algebra generated by a $n-$sequence. The sequences  $(p^-_0,\ldots,p^-_n)$ and $(p^+_0,\ldots,p^+_n)$ become monotonic and bounded stochastic processes with dependent increments.

\paragraph*{Uniformly sampled DOEs} The remainder of this article is devoted to a baseline statistical estimation of $\p$ in a monotonic framework, in a similar spirit to the MC approach in unconstrained frameworks. Therefore, in the following, the sampling is chosen uniform at each step: 
${\bf x_n}  \sim  {\cal{U}}_{\U_{n-1}}$.

\begin{figure}[hbtp]
\centering
     \input{basic6-bis.pstex_t}
\caption{Two-dimensional dominated and non-dominated subspaces after $n=14$ iterations. Points $\{0^2,\bf x_a,\bf x_b,\bf x_c,\bf x_d,\bf x_e,\bf x_f,\bf x_g\}$ have
nonzero signatures and are vertexes of $\U^-_n$. Points $\{\bf x_h,\bf x_i,\bf x_j,\bf x_k,\bf x_l,\bf x_m,\bf x_n,1^2\}$ have zero 
signatures and are vertexes of $\U^+_n$.} \label{graphic.example.2}

\vspace{1cm}

\input{basic11-quad.pstex_t}
\caption{Diagonal deterministic (DD-MRM) strategy, assuming a 
low $\p$, stopping after 4 steps.}
\label{types.of.surfaces}
\end{figure}

\section{A maximum likelihood estimator of $\p$}\label{MLE}

Assume that $\bf x_1,\ldots,\bf x_{n}$ are successively uniformly sampled in the nested non-dominated spaces $\U_0,\ldots,\U_{n-1}$. Next lemma follows.

\begin{lemma}\label{lemma.0} 
\hspace{0.25cm}\vspace{-0.5cm} ${\displaystyle p^-_n, p^+_n  \xrightarrow{a.s.}{}  \p.}$ \\
\end{lemma}

\noindent In corollary any estimator of $\p$ located between the bounds is strongly consistent. Especially, any crude average of the bounds gains a statistical validity. A more sophisticated approach can be carried out by noticing that, at step $k$, the occurence of a nonzero signature $\s_{\bf x_k}$ follows a Bernoulli distribution ${\cal{B}}(\gamma_k)$ conditionally to ${\cal{F}}_{k-1}$, with
\begin{eqnarray}
\nonumber \gamma_k & = & P\left(g({\bf x})\leq 0|{\bf x}\in\U_{k-1}\right), \\
\nonumber & = &  \frac{P\left(g({\bf x})\leq 0\right) - P\left(g({\bf x})\leq 0|{\bf x}\in\U^-_{k-1}\right)P\left({\bf x}\in\U^{-}_{k-1}\right)}{P\left({\bf x}\in\U_{k-1}\right)} 
\end{eqnarray}
from Bayes' formula, hence
\begin{eqnarray}
 \gamma_k & = & \frac{\p-p^-_{k-1}}{p^+_{k-1}-p^-_{k-1}}. \label{bernoulli}
\end{eqnarray}   
After $n$ steps, all information about $\p$ is  brought by 
the dependent-data likelihood $L_n(\p)=L_n(\p|{\bf x_1},\ldots,{\bf x_n})$ defined by the product of these  conditional Bernoulli
pdf: 
\begin{eqnarray}
L_n(\p) & = & \prod\limits_{k=1}^n
\left(\frac{\p-p^-_{k-1}}{p^+_{k-1}-p^-_{k-1}}
\right)^{\s_{\bf x_k}}\left(\frac{p^+_{k-1}-\p}{p^+_{k-1}-p^-_{k-1}}
\right)^{1-\s_{\bf x_k}}, \label{likelihood}
\end{eqnarray}
the maximum estimator (MLE) $\hat{\p}_n$ of which is considered in next proposition.  \\

\begin{proposition}\label{mle.exist}
Denote $\ell_n(\p)=\log L_n(\p)$. There exists a unique and consistent solution $\hat{\p}_n$ in $]p^-_{n-1},p^+_{n-1}[$ of the likelihood equation $\ell'_n(\p)  =  \sum_{k=1}^n \tilde{\omega}_k\left(\p\right) (p_k - \p)  = 0$, such that
\begin{eqnarray}
 \hat{\p}_n & = & \frac{\sum\limits_{k=1}^n
\tilde{\omega}_k\left(\hat{\p}_n\right) p_k}{\sum\limits_{k=1}^n
\tilde{\omega}_k\left(\hat{\p}_n\right)}, \label{mle} \\
\nonumber \text{with} \ \ \ \tilde{\omega}_k\left(\p\right) &  = & 
\left(\left(\p-p^-_{k-1}\right)\left(p^+_{k-1}-\p\right)\right)^{-1} \ \ \ \text{and} \ \ \   p_k  \ =  \ p^-_{k-1} + \left(p^+_{k-1}-p^-_{k-1}\right)\s_{\bf x_k} \label{unbiased.estimator.local.MLE}
\end{eqnarray}
\end{proposition}

Assumption 4 ensures the existence of $\hat{\p}_n$  since, by (\ref{cond001}) and (\ref{cond002}), $\p$ cannot be reached by at least one of the two bounds $(p^-_{n-1},p^-_{n+1})$ for any finite $n$. Similarly, the quantities defined in next propositions remain finite if the limit state surface $\Ss$ has mild smoothness properties. They are related to the behavior of the inverse of the Fisher information associated to (\ref{likelihood}), which converges to 0 faster than the variance of the usual MC $n-$estimator
\begin{eqnarray*}
V^{MC}_n(\p) & = & \frac{\p(1-\p)}{n}. 
\end{eqnarray*}

\begin{lemma}\label{lemma.expect} 
Assume that $\Ss$ is such that (\ref{cond001}) and (\ref{cond002}) hold ({\it Assumption 4}). Then, $\forall n\geq 0$,
\begin{eqnarray}
\E\left[1/(\p-p^-_{n})^2\right] & < & \infty, \label{expect.2} \\
\E\left[1/(p^+_{n}-\p)^2\right] & < & \infty, \label{expect.1} 
\end{eqnarray}
and consequently $\E[\tilde{\omega}_{n+1}(\p)]<\infty$.
\end{lemma}

\begin{proposition}\label{asympt.var.reduction}
Denote  
$J_n(\p)$ the Fisher information associated to (\ref{likelihood}). 
Then
\begin{eqnarray}
J^{-1}_n(\p) & = & \left(\sum\limits_{k=1}^n \E\left[\tilde{\omega}_k(\p)\right]\right)^{-1} \ \leq \  V^{MC}_n(\p) \frac{n}{\sum\limits_{k=1}^n (1-c_{k-1})^{-1}} \ < \ V^{MC}_n(\p) \label{exp2.var} \label{fisher.ecr}
\end{eqnarray}
where $c_0=0$ and $\forall~k>1$, 
\begin{eqnarray*}
c_k  &  = &  {\displaystyle \E\left[\frac{p^-_{k}}{\p}
+ \frac{1-p^+_{k}}{1-\p} -
\frac{p^-_{k}(1-p^+_{k})}{\p(1-\p)}\right].} \label{def.ck} \\
\nonumber
\end{eqnarray*}
\end{proposition}

\begin{proposition}\label{asympt.var.reduction.2}
Denote $\gamma_0=[(p^+_0-p^-_0)/p ^-_0]^2$. Then
\begin{eqnarray}
J^{-1}_n(\p) & \leq & V^{MC}_n(\p) \left(\frac{\p \gamma_0}{1-\p}\right). \label{exp1.var} \\
\nonumber
\end{eqnarray}
\end{proposition}

In this data-dependent context, the central limit Theorem \ref{asympt.norm.MLE} remains classical in the sense that the Cramer-Rao bound given by the inverse of the Fisher information is asymptotically reached by the MLE. It is technically based on the martingality of the score process $n\mapsto\{\ell'_{n}(\p)\}_n$.  Therefore inequalities (\ref{exp2.var}) and (\ref{exp1.var}) imply asymptotic  variance reduction with respect to Monte Carlo and robustness. From (\ref{exp1.var}), the asymptotic coefficient of variation (CV) of the MLE is such that
\begin{eqnarray*}
\mbox{CV}\left[\hat{\p}_n\right] & \leq &  \frac{\p}{\E\left[p^-_{n-1}\right]}\sqrt{\frac{\gamma_0}{n}} \ \overset{\infty}{\sim} \ \sqrt{\frac{\gamma_0}{n}}. \\ 
\end{eqnarray*}

\begin{theorem}\label{asympt.norm.MLE} 
Let $(\lambda_n)$ be any deterministic sequence in $]0,1[$ such that $\lambda_n \to 1$. Under the supplementary assumptions:
\begin{description}
\item[(i)] 
${\displaystyle \frac{1}{n^{\delta}}\sum\limits_{k=1}^n \left(\tilde{\omega}_k(\p) - \E\left[\tilde{\omega}_k(\p)\right] \right)  \xrightarrow[]{\Pp}   0} \ \ \ \text{for any $\delta\geq 1$.}$ 
\item[(ii)] ${\displaystyle \frac{p^+_n-\p}{\p - p^-_n} \xrightarrow{\Pp}{} 1.}$
\item[(iii)] ${\displaystyle \frac{\bar{p}_n - \p}{p^+_n - \p} \xrightarrow[]{\Pp}   0 \ \ \  \text{and} \ \ \  \frac{\bar{p}_n - \p}{\p - p^-_n} \xrightarrow[]{\Pp}   0} \ \ \ \text{with $\bar{p}_n=(1-\lambda_n)\hat{\p}_n + \lambda_n\p$}$
\end{description}
then 
\begin{eqnarray}
\hspace{2cm} {\displaystyle {{J^{1/2}_n(\p)}} \left(\hat{\p}_n - \p\right)  \xrightarrow[]{{\cal{L}}}  {\cal{N}}(0,1).} \label{tlc1} 
\end{eqnarray}  
\end{theorem}
\vspace{0.15cm}

The law of large numbers {\bf (i)} reflects the requirement that the sum of weights $\tilde{\omega}_k(\p)$ cannot diverge faster than $\mathcal{O}(J_n(\p))$ from its mean behavior when $n\rightarrow\infty$.  Although difficult to check in practice, this behavior seems rather intuitive because the sampling trajectories mainly vary at the first steps of the algorithm, when the non-dominated space is still large. 
 Therefore {\bf (i)} can be perceived as an indirect requirement on the surface $\Ss$. 
 Assumption {\bf (ii)} appears somewhat natural, saying that the bounds converge to $\p$ symmetrically. Assumption {\bf (iii)} expresses the idea that any estimator located between $\hat{\p}_n$ and $\p$ converges to $\p$ faster than the bounds. Again, it seems intuitive since $\hat{\p}_n$ is defined as an incremental average (cf. (\ref{mle})), and therefore adopts a smoother behavior than the bounds, as a function of $n$. 

Next proposition allows for an empirical estimation of the asymptotic variance and confidence intervals. The additional requirement {\bf (v)} appears mild and in the same spirit than the smoothness assumptions on $\Ss$, saying that $\p$ cannot be exactly reached by an average of the bounds for any finite number $n$ of trials.

\begin{proposition}\label{TCL.variance}
Denote $\hat{J}_n(\p)=\sum_{k=1}^n \tilde{\omega}_k(\p)$. Under the assumptions of Theorem \ref{asympt.norm.MLE}, and assuming in addition:
\begin{description}
\item[(iv)] Assumption {\bf (i)} remains true $\forall \delta\geq 1/2$,
\item[(v)]  $\nexists~n<\infty$ such that 
$\p  =  \left.(2n)^{-1} \sum_{k=1}^n \tilde{\omega}_k(\p) (p^-_{k-1} + p^-_{k-1})\right/\sum_{k=1}^n \tilde{\omega}_k(\p)$,
\end{description}
then
\begin{eqnarray}
\hspace{2cm} \frac{\hat{J}_n^{~5/2}(\p)}{|\hat{J}'_n(\p)|}\left(\hat{J}_n^{~-1}(\hat{\p}_n)-J^{-1}_n(\p)\right) & \xrightarrow[]{{\cal{L}}} & {\cal{N}}(0,1). \label{tlc2} 
\end{eqnarray}
\end{proposition}
\vspace{0.15cm}

The reality of the theoretical descriptions hereinbefore is examined in the two next sections, through numerical experiments conducted on toy examples and a more realistic hydraulic model.

\section{Numerical experiments I: toy examples}\label{numerical.experiments}

The statistical behavior of the MLE is illustrated here using the following 
 generic toy example. For a given dimension $d$, denote
\begin{eqnarray*}
Z_d & = & h_d({\bf Y}) \ = \ {Y_1}/{(Y_1 + \sum\limits_{i=2}^d Y_i)}
\end{eqnarray*} 
where the physical input $Y_i$ follows the gamma distribution ${\cal{G}}(i+1,1)$ with cdf $F_{Y_i}$, independently of other inputs. Obviously, $\forall~d\geq 2$, $h_d$
is increasing in $(-X_1,X_2\ldots,X_d)$ where $X_i=F_{Y_i}(Y_i)\sim{\cal{U}}_{[0,1]}$, and 
$Z_d$ follows the beta distribution ${\cal{B}}_e(2,2^{-1}{(d+1)(d+2)} - 3)$. Therefore, denoting
$q_{d,\p}$ the $\p-$order quantile of $Y_d$, the deterministic function defined by
\begin{eqnarray*}
g_d({\bf X}) & = h_d\circ T^{-1}({\bf X}) - q_{d,\p},
\end{eqnarray*}
with $T^{-1}({\bf x})=(F^{-1}_{Y_1}(x_1),\ldots,F^{-1}_{Y_d}(x_d))$, 
is related to the known exceedance probability $\p$. 

\Subsection{First results}

\begin{figure}[hbtp]
\centering
\begin{minipage}[c]{.46\linewidth}
      \includegraphics[width=7cm,height=7cm]{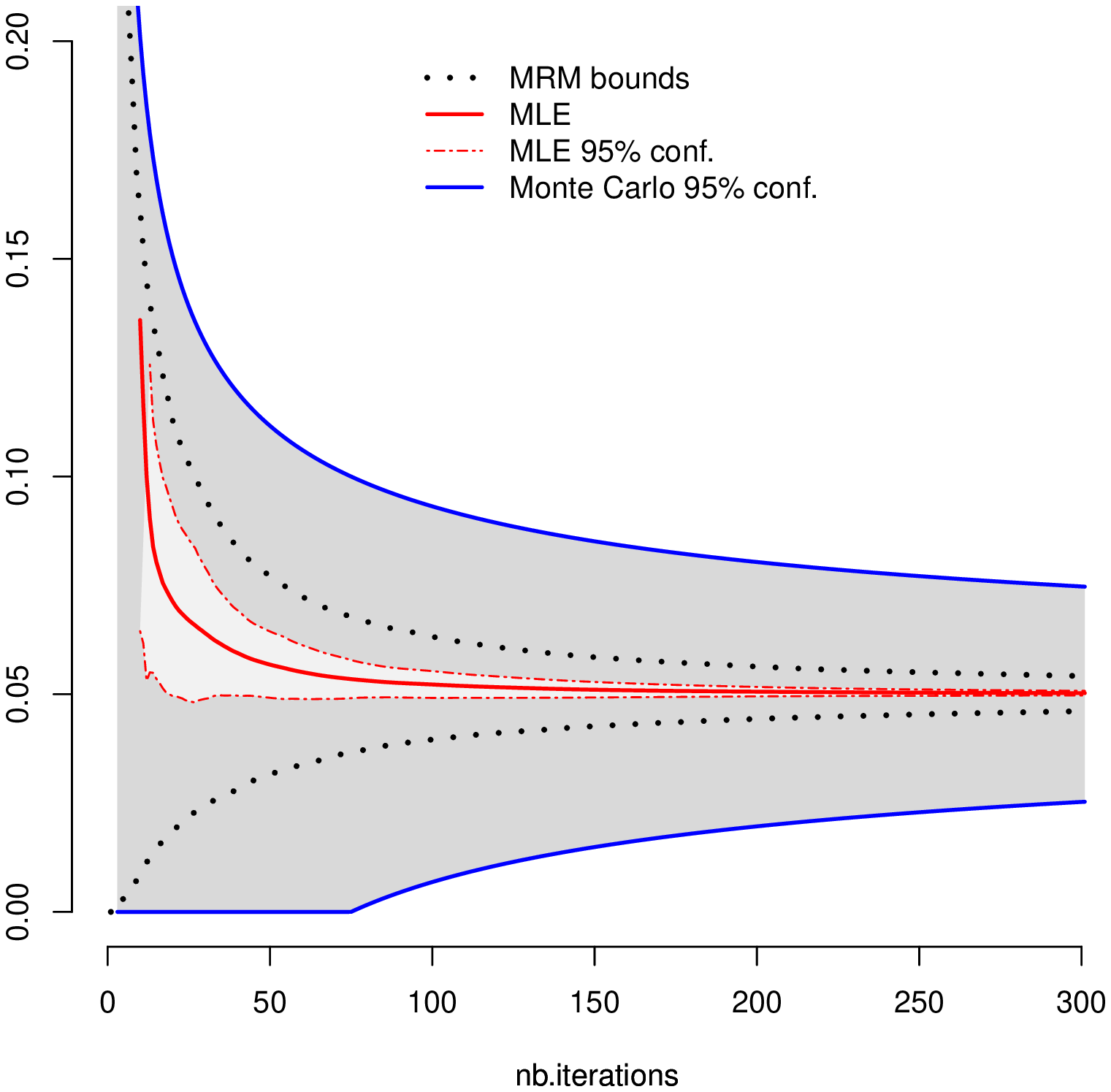}
\caption{MRM deterministic bounds and MLE, with Monte Carlo and MLE 95\%-confidence areas, in dimension $d=2$, for $\p=5\%$.  Empirical estimations are made  over 300 parallel MRM trajectories.}
\label{bounds-MLE}
   \end{minipage} \hfill
   \begin{minipage}[c]{.46\linewidth}
   \vspace{-1.25cm}
      \includegraphics[width=7cm,height=7cm]{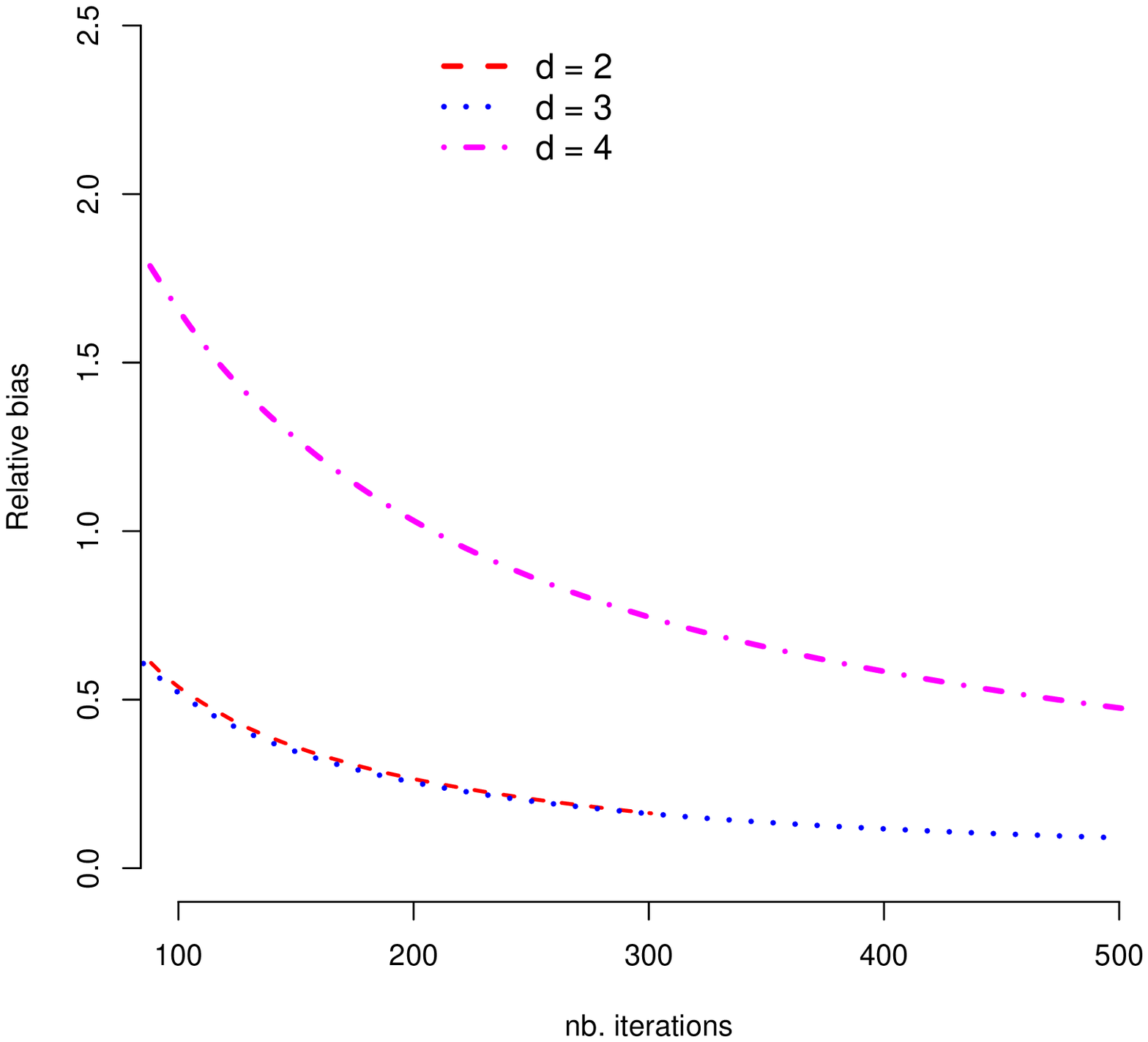}
\caption{Relative bias of the MLE $\hat{p}_n$ for the dimensions $d\in\{2,3,4\}$.}
\label{relative-bias-MLE}
   \end{minipage}

\vspace{1cm}

\includegraphics[width=14cm,height=7cm]{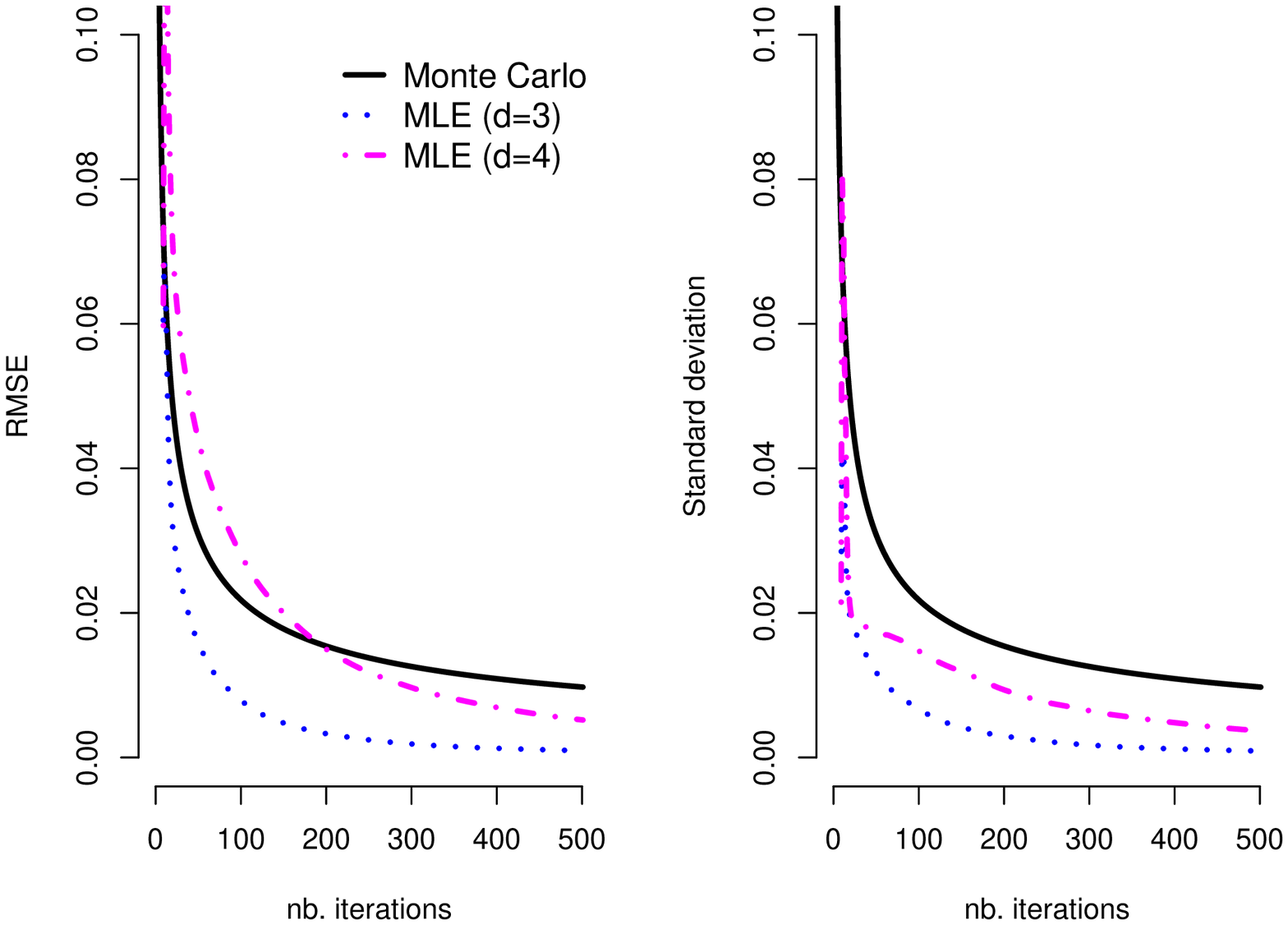}
\caption{Root mean square error (RMSE; {\it left}) and standard deviation ({\it right}) of the standard Monte Carlo estimator and $\hat{\p}_n$ for $d=3$ and $d=4$, for $\p=0.05$. Empirical estimations are made  over 300 parallel MRM trajectories.}
\label{RMSE-var-MLE}
\end{figure}

\begin{figure}[hbtp]
\centering
\includegraphics[width=7cm,height=7cm]{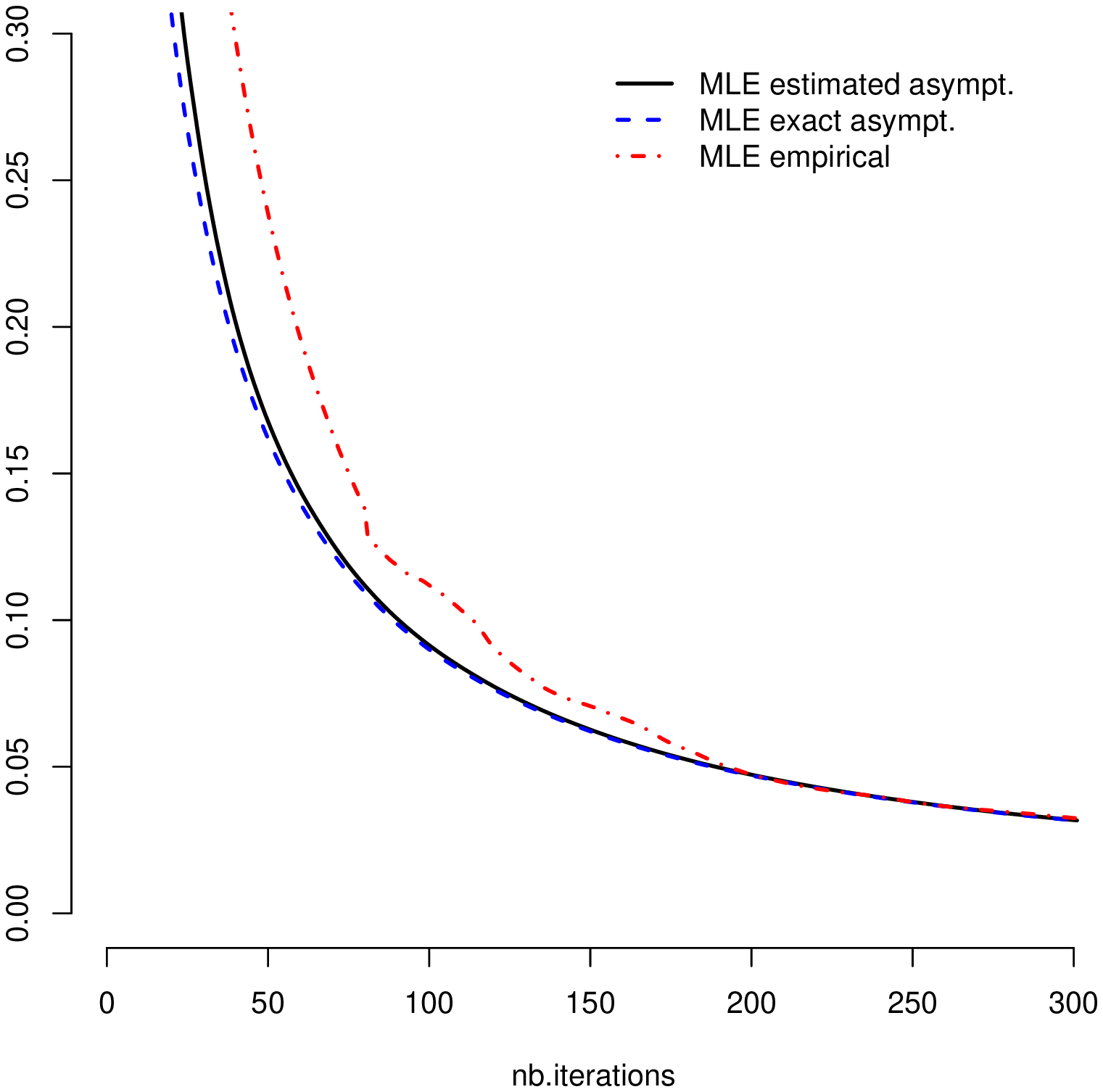}
\caption{Ratios of MLE standard deviations over Monte Carlo standard deviations, computed over 100 MRM replications, in dimension $d=2$.}
\label{ratios}

\vspace{1cm}

   \input{basic14.pstex_t}
   \vspace{-1.5cm}
\caption{Two-dimensional situation after $n=14$ iterations. A replication $\hat{\Ss}_n$ (dashed curve) of $\Ss$ can be produced based on a monotonic neural network prediction of $\s_{\bf x}$ in $\U_n$. The volume under $\hat{\Ss}_n$ corresponds to a new probability $\tilde{\p}_n$ in the magnitude of $\p$, which can be estimated by Monte Carlo at an arbitrary precision.} 
\label{emulation.fig}
\end{figure}

In dimension 2, using $\p=5\%$,  the behavior of MRM bounds can be easily compared to the MC 95\%-confidence area (Figure \ref{bounds-MLE}). This small dimension induces a significant improvement in precision with respect to Monte Carlo, which however disappears in higher dimensions and highlights the need for real statistical estimators. Studies of the root mean square error (RMSE) and the standard deviation of the MLE, which are plotted in Figure \ref{RMSE-var-MLE} as functions of the increasing number of calls to $g_d$ for dimensions 3 and 4, reflected the high variance reduction of the iterative estimator $\hat{p}_n$ with respect to Monte Carlo but highlighted a positive bias (Figure \ref{relative-bias-MLE}). Indeed the highest weights favor local estimators $p_k=p^+_k$ when approaching  $\Ss$ (ie., when $\s_{\bf x_k}=1$ in (\ref{unbiased.estimator.local.MLE})).  On the examples considered in this last figure (as well as in other experiments not shown here), a marked gap in relative bias was noticed between dimensions 3 and 4. Under dimension 4, the bias remains reasonable from a moderate number of iterations (typically 400). Else it dramatically stays at high values. Other experiments have shown on this example the effective convergence of the empirical variance of the MLE towards the Cramer-Rao bound as well as the good behavior of its empirical estimate (Figure \ref{ratios}).

\Subsection{Bias correction via bootstrap heuristics}

Bias removal appears as a practical requirement, automatizing the estimation of $\p$. Indeed,  given a finite value of $n$, estimating $\p$ requires to decide from which iteration $k_n\geq 1$ the computation of the MLE $\hat{\p}_n$ can be worth it, redefining
$\hat{\p}_n  =  \sum_{i=k_n}^n \psi_{i,n}(\hat{\p}_n) p_i$ 
with $\psi_{i,n}(\p)=\tilde{\omega}_i(\p)/\sum_{j=k_n}^n \tilde{\omega}_j(\p)$. An intuitive rule is to select
\begin{eqnarray}
k^*_n & = & \arg\min\limits_{k_n} \mbox{RMSE}\left(\hat{\p}_n\right).\label{selecting}
\end{eqnarray} 
If the MLE were debiased, $\mbox{RMSE}\left(\hat{\p}_n\right)\simeq J^{-1}_n(\p)$  which  is minimized by $k^*_n=1$. 

Given a fixed number $n$ of trials, two general approaches may be used for controlling and correcting the bias ${B}_{n}=\E[\hat{\p}_n]-\p$ affecting $\hat{p}_n$. A {\it corrective} approach consists of obtaining a closed-form expression for the bias from Taylor expansions \citep{COX1968,FER98} or penalizing the score or the likelihood functions \citep{BESTER2005}. This approach is not carried out here since the data-dependent context would require a specific algebraic work and technical developments about the empirical estimation of the main quantities involved. The alternative use of bootstrap resampling techniques \citep{Efron1993} can assess the bias empirically. 
For the simplicity of their principle, these heuristics are preferred here.  
 
 In the present context, bootstrap experiments must be  based on a replication $\hat{\Ss}_n$ of the limit state surface $\Ss$ (see Figure \ref{emulation.fig} for an illustration). Under Assumption 4, $\Ss$ can be interpreted as the decision frontier of a supervised classification binary problem, without horseriding of classes (ie., perfectly separable). Therefore $\hat{\Ss}_n$ depends on the choice of a classifier $\hat{C}_{n,M}$ calibrated from an arbitrary number $M$ of points sampled in dominated subspaces. The rationale of the bootstrap heuristics is as follows. Given $\hat{C}_{n,M}$, 
 the signature $\s_{\bf x}$ of any ${\bf x}\in\U_n$ can be predicted by the occurence of 	
$P(g({\bf x}\leq 0)|\hat{C}_{n,M})\geq 1/2$. Then denote $\tilde{\p}_{n,M}$ the volume under $\hat{\Ss}_{n}$. It can easily be estimated by $\tilde{p}_{n,M,Q}$ at an arbitrary precision by Monte Carlo sampling (depending on $Q$). Moreover a large number $S$ of MLE estimators of $\tilde{\p}_{n,M}$ can be fastly computed using the predicted signatures. \\

The bootstrap heuristics make sense if the classifier is chosen such that $\tilde{\p}_{n,M}\to\p$ when $(n,M)\to\infty$, so that the features of the experiment are asymptotically reproduced. Moreover, $\hat{C}_{n,M}$ must produce a monotonic surface $\hat{\Ss}_n$. For these reasons, the four-layer monotonic Multi-Layer neural networks (MLNN) proposed by \citet{DAN10} have been chosen for the experiments. Based on a combination of minimum and maximum functions (so-called MIN-MAX networks) over the two hidden layers, these networks have universal approximation capabilities of  monotonic continuous functions. Besides, this choice matches the advices by \citet{HUR04} who strongly recommended the MLNN and Support Vector Machines (SVM) to estimate $\Ss$ in a structural reliability context. Both tools are flexible, can estimate a frontier on the basis of a few samples and overcome the curse of dimensionality. \\

\clearpage
\texttt{\hspace{2.5cm}{\bf Classification-based bootstrap algorithm}
\begin{enumerate}
\hbox{\raisebox{0.4em}{\vrule depth 0pt height 0.4pt width 13.5cm} }
\item Sample  $\bf x^+=(\bf x^+_1,\ldots,\bf x^+_M) \overset{iid}{\sim} {\cal{U}}_{\U^+_n}$ and $\bf x^+=(\bf x^-_1,\ldots,\bf x^-_M)\overset{iid}{\sim} {\cal{U}}_{\U^-_n}$. 
\item From $(\bf x^+,\bf x^-)$, build a monotonic classifier $\hat{C}_{n,M}$ of $(\U^-,\U^+)$. 
\item Replace $g$ by the uncostly monotonic (increasing) function
\begin{eqnarray*}
\hspace{2cm} \tilde{g}({\bf x}) & = & \left\{ \begin{array}{ll} -1 & \text{if $P(g({\bf x}\leq 0)|\hat{C}_{n,M})\geq 1/2$,} \\ +1 & \text{else.} \end{array}\right. 
\end{eqnarray*}
\vspace{-0.2cm}
\item Sample $\bf x_{1},\ldots,\bf x_{Q} \overset{iid}{\sim} {\cal{U}}_{\U}$ and compute 
 $\tilde{p}_{n,M,Q}  =  Q^{-1}\sum_{k=1}^Q \1_{\{\tilde{g}({\bf x_k})\leq 0\}}$.
\item For $i=1,\ldots,S,$ get a MLE estimator $\tilde{p}^{(i)}_{n,M,Q}$ then estimate ${B}_{n}$ by
\begin{eqnarray*}
\hspace{2cm} \hat{B}_{n,M,Q,S}  & =  & S^{-1}\sum_{i=1}^S \tilde{p}^{(i)}_{n,M,Q} - \tilde{p}_{n,M,Q}.
\end{eqnarray*}
\end{enumerate}
\vspace{-0.2cm}
\hbox{\raisebox{0.4em}{\vrule depth 0pt height 0.4pt width 13.5cm} }
}
\vspace{0.25cm}

Numerical tests in function of $n$ and $d$ were conducted, and the results are presented in Table \ref{experimental.results.1}. The bias correction is found to be effective even from a moderate number of iterations (some hundreds) until dimension 5, and a budget of at least $n=1,000$ is enough to correct a bias in dimension 8. With less than 10\% of overestimation on average on this example, these bootstrap heuristics also appear relevant when the exact value of $\p$ is less interesting than its magnitude, which is often the case in design optimization where it is aimed to diminish $\p$ of a given factor by constraining the inputs \citep{Tsompanakis2007}.

\begin{table}[hbtp]
\centering
\begin{tabular}{cccccccccccc}
\hline
& \multicolumn{10}{c}{Dimension $d$} & \\
& &&&&&&&&&&\\
\cline{3-12} 
   &  \multicolumn{5}{c}{$\p=0.05$} & &  \multicolumn{4}{c}{$\ \ \ \p=0.005$} & \\
   & &&&&&&&&&&\\
$n$   &&  2 & 3 & 4 & 5 & 8 &&  2 & 3 & 4 & 5 \\
\cline{3-7}\cline{9-12} 
& &&&&&&&&&&\\
\small 50   && \small 1.19  & \small 4.17 & \small  6.93 & \small 16.87 & 27.85  && \small 8.21  & \small 13.80   & \small 18.22 & \small 39.57 \\ 
\small 100  && \small 0.28  & \small 2.31 & \small  4.79 & \small 12.94 & \small 22.98  && \small 6.15  & \small 11.55   & \small 15.67 & \small 31.40  \\
\small 250  && \small 0.21  & \small 1.87 & \small  3.34 & \small 8.74  & \small 19.12  && \small 3.28  & \small 8.72    & \small 11.01 & \small 24.68  \\
\small 500  && \small 0.12  & \small 1.25 & \small  2.87 & \small 6.20  & \small 16.76  && \small 1.14  & \small 5.84    & \small 8.12  & \small 16.06  \\
\small 1000 && \small -0.02  & \small 0.47 & \small 2.14 & \small 2.97  & \small 12.85  && \small 0.12 & \small 2.72    & \small 5.28  & \small 9.23  \\
\small 2000 && \small -0.008 & \small -0.28 & \small  1.61 & \small 2.08  & \small 7.66   && \small -0.34 & \small 1.55  & \small 3.09  & \small 6.51 \\
\hline
\end{tabular}
\vspace{0.5cm}
\caption{Relative error in \% between estimated bias and real bias, for two true probabilities $p=5\%$ and $p=0.5\%$. Results are averaged on 100 experiments, each boostrap estimation being based on $S=1,000$ MLE replicates. For each $n$, the neural network is build from $M=10^6$ sampled vectors, with a classification error rate less than 0.25\% on these training data.}
\label{experimental.results.1}
\end{table}

\section{Numerical experiments II: a simplified hydraulic case-study}

De Rocquigny (2009), \citet{LIM10} then \citet{MUN10} considered a simplified but realistic hydraulic model linking the downstream water level $H$ $(m)$ of a river section, of width $b=300$ $(m)$ and length $l=5000$ $(m)$, with the upstream discharge $Q$ $(m^3/s)$ and the friction  coefficient $K_s$ $(m^{1/3}/s)$ of the river bed. Denoting $Z_m$ and $Z_v$ the upstream and downstream altitude of the river bed above seal level,
\begin{eqnarray*}
H & = & \left(\frac{Q}{bK_s\sqrt{\frac{Z_m-Zv}{l}}}\right)^{3/5}.
\end{eqnarray*}
Assuming a dike level $h_0=55.5$ (m), the flood probability is $\p=P(g'\circ T^{-1}({\bf X})\leq 0)$ where ${\bf Y}=\{Q,K_s\}$ and $T=(F_Q,F_{K_s})$ (2-dim. version) or ${\bf Y}=\{Q,K_s,Z_m,Z_v\}$ and $T=(F_Q,F_{K_s},F_{Z_m},F_{Z_v})$ (4-dim. version),  and
\begin{eqnarray*}
g'({\bf Y}) & = & h_0 - Z_v - H({\bf Y}),
\end{eqnarray*} 
which is increasing in $(-Q,K_s,Z_m,-Z_v)$. Input distributions or punctual values are chosen as in \citet{LIM10}. $Q$ follows a Gumbel distribution with location 1013 and scale 558, truncated in $[10,10^4]$. $K_s$ is normal ${\cal{N}}(27.8,3^2)$ truncated in 0. In the 4-dim. version, $Z_m$ and $Z_v$ are triangular on $[53.5,56.5]$ and $[48.5,51.5]$ with respective modes 55 and 50 (their respective values in the 2-dim. version). \\

For several computational budgets and averaged over 100 repeated experiments, two alternative methods are compared to the MRM bias-corrected MLE: the MC method and an engineering FORM-IS method build on two steps: {\bf (a)}  with a limited number of trials (no more than 40), the First-Order Reliability Method (FORM) is run to provide an estimate of the conception point $\bm\beta=\arg\min\|{\bf u}\|$ on $\{g'\circ F^{-1}\circ\Phi({\bf u})\leq 0\}$, with $\Phi$ the standard normal pdf and $\bf u$ a random variable evolving in the $d-$dimensional standard Gaussian space $U$ ; {\bf (b)} an Importance Sampling (IS) method that uses the budget left to sample in $U$ using a standard normal distribution centered on $\bm\beta$. See \citet{OPENTURNS2011} for details about the implementation of the method.  

For a given $n$, the three methods are compared through the following indicators: $\E[\hat{\p}_n]$, $\mbox{CV}[\hat{\p}_n]$ and the relative average precision $\gamma_n=\E[(p^+_n-p^-_n)]/\p$. Using the DOEs produced by the MC and the FORM-IS methods, these bounds can obviously be computed accounting for the monotonicity of $g$. For each version a MC computation involving 40,000 particules provides a precise estimate of $\p$, which is used for estimating $\p$ in $\gamma_n$. Finally, $S=1,000$ bootstrap replicates were used for the correction of each MRM-MLE estimate. The results are summarized on Table \ref{24dim-tab}. \\

In terms of magnitude, the three methods perform similarly. The benefit of using MRM instead of MC or FORM-IS in these low dimensions clearly appears in most cases, and more obviously in dimension 2: MC needs at least 200 times more iterations than MRM to reach a similar precision $\mbox{CV}[\hat{\p}_n]$, and if FORM-IS is significantly better than MC, the precision of its estimates remains far beyond of those produced by MRM.  In dimension 4, the difference between these two methods somewhat vanishes and they lead to close performance when the number of iterations remains low. For both dimensional cases, it was  noticed that a single FORM run can provide a crude estimate of $\p$ with good magnitude after 10 iterations only. But the dimensional increasing allows the part of the importance sampling falling into the non-dominated area to be greater than in a two-dimensional setting.  \\  

\begin{table}[hbtp]
\centering
\begin{tabular}{cl|lllclll}
\hline
$n$ & method & \multicolumn{3}{c}{dimension = 2} && \multicolumn{3}{c}{dimension = 4} \\  
    &        & \small $\E[\hat{\p}_n]$ & \small $\mbox{CV}[\hat{\p}_n]$  & \small $\gamma_n$  &&  \small $\E[\hat{\p}_n]$ & \small $\mbox{CV}[\hat{\p}_n]$  & \small $\gamma_n$  \\ 
\cline{3-5} \cline{7-9}
& \\
\small 100  & \small MC     &  \small 0.002775 &\small 190\% &\small  2,900\% &&\small 0.010075 &\small 99\% & \\
    & \small FORM-IS        & \small 0.002241 &\small 68\% &\small 478\%     &&\small 0.018147 &\small 74\% & \\
    & \small MRM            & \small 0.002781 &\small 14\% &\small 48\%      &&\small 0.015498 &\small 82\%  &\small 1,400\%\\
& \\
\small 200 & \small MC      &\small 0.002775 &\small 134\% &\small 630\%     &&\small 0.010075 &\small 70\% &\small 2,300\%\\
    & \small FORM-IS        &\small 0.002667 &\small 44\%  &\small 244\%     &&\small 0.010242 &\small 42\% &\small 2,230\%\\
    & \small MRM            &\small 0.002776 &\small 6\%   &\small 24\%      &&\small 0.012451 &\small 35\% &\small 800\% \\  
& \\
\small 1,000 & \small MC    & \small 0.002775 &\small 60\% &\small 515\%     &&\small 0.010075 &\small 31\% & \small 1,200\% \\
    & \small FORM-IS        & \small 0.002736 &\small 27\% &\small 168\%     &&\small 0.009959 &\small 27\% & \small 1,000\%\\
    & \small MRM            & \small 0.002775 &\small 0.12\% &\small 5.6\%   &&\small 0.010911 &\small 20\% & \small 300\% \\ 
& \\
\small 40,000 & \small MC      &\small  0.002775 &\small 9.5\% &\small 475\%  &&\small 0.010075 &\small 5\% &\small 247\% \\
\hline
\end{tabular}
\vspace{0.5cm}
\caption{Estimation results for the two-dimensional  and four-dimensional versions of the problem.}
\label{24dim-tab}
\end{table}


\section{Discussion}

Many structural reliability problems deal with the fast estimation of a probability $\p$ of an undesirable event. This event can often be defined by the occurence of an exceedance  in 
output of some time-consuming function $g$ with stochastic multidimensional inputs. In the present article, $g$ is assumed to be monotonic and possibly non-continuous.  

Pursuing pioneering works by \citet{ROC09} and \citet{LIM10} who explored heuristically the benefits of this framework, this article  first offers a formal description of the latter that focuses on the existence of deterministic bounds around $\p$. A sequential strategy of numerical experiments in the input space allows for a progressive narrowing of this interval. The second and main aspect of the paper is the definition and the study of a statistical estimator of $\p$ when the strategy becomes stochastic and leans on uniform nested sampling.  Easy to compute, it is defined as the maximizer of a likelihood (MLE) of dependent data sampled from Bernoulli distributions, whose parameters are explicit functions of the dynamic bounds.  

A keypoint of the paper is the theoretical description of its asymptotic properties, which are found similar to those arising from the classical estimation theory, provided some intuitive assumptions are respected. They are found mild in practice on some examples. Both theoretical and applied results show a significant improvement of the fastness and the robustness of this estimator with respect to the usual Monte Carlo estimator. In the third part of the paper, boostrap heuristics are proposed and carried out successfully to remove the non-asymptotic bias affecting the MLE, via constrained neural networks. Only a basic continuity assumption on the limit state (or {\it failure}) surface is needed to benefit from their universal approximation capabilities. \\

Thus, the tools proposed in this article and its supplementary material in Appendix can be directly used in structural reliability applications, without preliminary learning step (as usual, for instance, in stratified methods). However, the generality of the frame allows for a wider range of theoretical and applied studies. These research avenues are briefly discussed in the following items. \\

\paragraph*{Bias correction} Following \citet{HUR04}, support vector machines (SVM) should probably be considered instead of neural networks, since their geometric interpretation of margin maximizers appears more suitable.  In addition to the monotonicity constraint, they should be build at step $n$ under the linear  constraint that the volume under the predicted surface be equal to the current (biased) estimator $\hat{\p}_n$. This would certainly improve the properties of the bootstrap heuristics. More importantly, this method should be now tested on a large variety of examples, and the intuitive feeling of its ability to    correct the bias must be confirmed by more applied and theoretical studies. 

In parallel, future studies should focus on adopting a corrective approach to the bias affecting the MLE, then on selecting a slippery window of indexes, according to (\ref{selecting}) or a similar rule, such that the MLE converges faster to $\p$. The comparison of the experimental benefits of both approaches would help the method to become more ready-to-use.      

\paragraph*{Simplifying the assumptions} Most of the technical assumptions that are needed to get the theoretical results present some intuitive features, and are underlyingly linked to the nature of the limit state surface. However, they remain difficult to check in practice, although asymptotic normality was always noticed in  numerical experiments. Therefore, future work should be dedicated to simplifying those assumptions and classifying the limit state surfaces in function of their ability to allow a fast and robust estimation of $\p$. 

\paragraph*{Sensitivity studies} Crucial tasks in structural reliability are sensitivity studies of probabilistic indicators to the uncertainty input model \citep{MORIO2011}. Therefore, assuming ${\bf X}=T({\bf Y})$ where the $\bf Y$ represent {\it physical} inputs with multivariate distributional transform $T({\bf y}=(F_1(y_1),\ldots,F_d(y_d))$ (each $F_i$ being the marginal cdf of $Y_i$), given a budget $n$, the variations of $(p^-_n,p^+_n,\hat{\p}_n)$ due to modifying $T$ in $T_{\epsilon}$ should be the subject of future works. As a supplementary benefit of the method, the new values $(p^-_{k,\epsilon},p^+_{k,\epsilon},\hat{\p}_{n,\epsilon})$, for $k\in\{0,\ldots,n\}$, can be recomputed without any supplementary call to $g$, thanks to an importance sampling mechanism. Indeed, as the subspaces $(\U^-_k,\U^+_k)$ remain dominated whatever the choice made on input distributions in the physical space, then
\begin{eqnarray*}
p^-_{k,\epsilon} & = & {\displaystyle \int_{T^{-1}(\U^-_k)} d T_{\epsilon}({\bf y})} \ \ \ \text{and} \ \ \ 
p^+_{k\epsilon} \ = \ {1-\displaystyle \int_{T^{-1}(\U^+_k)} d T_{\epsilon}({\bf y})}, 
\end{eqnarray*}
which can computed by a simple Monte Carlo method. In such future studies, we suggest that the progressive bounds could be defined as {\it robust} if they remain true whatever the fluctuations of $F_{\epsilon}$ in a well-funded variational class around $F$.   

\paragraph*{Exploring other forms of stochastic DOEs} A keypoint of future works will be to elaborate unbiased estimators from sequential stochastic designs of experiments with non-asymptotic properties. Indeed, the asymptotic variance of the MLE reaches the Cramer-Rao bound $J^{-1}_n(\p)$. Therefore any unbiased estimator based on sequential uniform sampling, especially those defined by
$\tilde{\p}_n  =  \sum_{k=1}^n \omega_k p_k$ where the $\omega_k$ are now {\it deterministic} weights, independent on $\p$ and summing to $1$, will never reach a lower variance than 
$J^{-1}_n(\p)$, even though the $\omega_k$ are optimized. Improving the Monte Carlo acceleration $nJ^{-1}_n(\p)/(\p(1-\p))$ will only be possible using less naive strategies than uniform samplings. The problem of defining such samplings so that an unbiased estimator of $\p$ has better statistical properties will be the subject of a future paper. 

\paragraph*{Towards partial monotonicity} Finally, the practical limits of monotonicity assumptions should be refined.  Intuitively, monotonicity as a building hypothesis seems  antagonist to high-dimensional structural safety problems, and could mainly characterizes the behavior of $g$ as a function of its most influential input variables (as determined by global sensitivity analyses). Indeed, the real examples treated by \citet{ROC09,LIM10} and \citet{RAJ10} do not go beyond dimension 4.   
Partial monotonicity, as defined by \citet{DAN10}, is a more appealing and realistic property, for which the methods developed in a pure monotonicity context should be adapted in the future.


\section{Acknowledgements}

 E. de Rocquigny (\'Ecole Centrale de Paris) must be thanked to have drawn my attention to this topic. I am especially grateful to all members of the MRI/T56 group and A. Dutfoy-Lebrun at EDF Research \& Development and Prof. Fabrice Gamboa and Dr. Thierry Klein (Institut de Math\'ematiques de Toulouse), for their support during this work. Finally, I am thankful to an anonymous reviewer who deeply read this article and made comments, critics and useful suggestions which greatly helped to improve it.

%
%
%
%
\bibliography{bibliotheque}

\begin{thebibliography}{39}
\providecommand{\natexlab}[1]{#1}
\providecommand{\url}[1]{\texttt{#1}}
\providecommand{\urlprefix}{URL }

\bibitem[{Anonymous(2011)}]{OPENTURNS2011}
Anonymous (2011) Open {TURNS} version 0.14.0 - {R}eference guide. Tech. rep.,
  EDF - EADS - PhiMeca.

\bibitem[{Bercu(2008)}]{BER08}
Bercu, B. (2008) Inégalités exponentielles pour les martingales [{\it in
  {f}rench}]. \emph{Journées ALEA 2008 CIRM} pp. 10--14 March.

\bibitem[{Bester \& Hansen(2005)}]{BESTER2005}
Bester, C. \& Hansen, C. (2005) Bias reduction for {B}ayesian and frequentist
  estimators. {W}orking {P}aper. \emph{University of Chicago} .

\bibitem[{Cannamela \emph{et~al.}(2008)Cannamela, Garnier \& Iooss}]{CAN08}
Cannamela, C., Garnier, J. \& Iooss, B. (2008) Controlled stratification for
  quantile estimation. \emph{Annals of Applied Statistics} \textbf{2},
  1554--1580.

\bibitem[{Chan(2008)}]{CHA08}
Chan, T. (2008) A (slightly) faster algorithm for {K}lee's measure problem. pp.
  94--100, College Park, MD, USA.

\bibitem[{Chen(2009)}]{CHEN09}
Chen, G. (2009) Monotonicity of dependence concepts: from independent random
  vector into dependent random vector. \emph{World Academy of Science,
  Engineering and Technology} \textbf{57}, 399--408.

\bibitem[{Chlebus(1998)}]{CHL98}
Chlebus, B. (1998) On the {K}lee's measure problem in small dimensions.
  \emph{Proceedings of the 25th Conference on Current Trends in Theory and
  Practice of Informatics} .

\bibitem[{Cox \& Snell(1968)}]{COX1968}
Cox, D. \& Snell, E. (1968) A general definition of residuals. \emph{Journal of
  the Royal Statistical Society} \textbf{30}, 248--275.

\bibitem[{Crowder(1975)}]{CRO75}
Crowder, M. (1975) Maximum likelihood estimation for dependent observations.
  \emph{Journal of the Royal Statistical Society} \textbf{38}, 43--53.

\bibitem[{Crowder(1983)}]{CRO83}
Crowder, M. (1983) On constrained maximum likelihood estimation with non-iid.
  observations. \emph{Annals of the Institute of Statistical Mathematics}
  \textbf{36}, 239--249.

\bibitem[{Daniels \& Velikova(2010)}]{DAN10}
Daniels, H. \& Velikova, M. (2010) Monotone and partially monotone neural
  networks. \emph{IEEE Transactions in Neural Networks} \textbf{21}, 906--917.

\bibitem[{de~Berg \emph{et~al.}(1997)de~Berg, van Kreveld, Overmars \&
  Schwarzkopf}]{deB97}
de~Berg, M., van Kreveld, M., Overmars, M. \& Schwarzkopf, O. (1997)
  \emph{Computational {G}eometry {A}lgorithms and {A}pplications}.
  Springer-Verlag.

\bibitem[{de~Rocquigny(2009)}]{ROC09}
de~Rocquigny, E. (2009) Structural reliability under monotony: A review of
  properties of {FORM} and associated simulation methods and a new class of
  monotonous reliability methods ({MRM}). \emph{Structural Safety} \textbf{31},
  363--374.

\bibitem[{Durot(2008)}]{DUROT08}
Durot, C. (2008) Monotone nonparametric regression with random design.
  \emph{Mathematical Methods in Statistics} \textbf{17}, 327--341.

\bibitem[{Efron \& Tibshirani(1993)}]{Efron1993}
Efron, B. \& Tibshirani, R.J. (1993) \emph{An {I}ntroduction to the
  {B}ootstrap}. New York: Chapman \& Hall.

\bibitem[{Erickson(1998)}]{ERI98}
Erickson, J. (1998) {K}lee's measure problem. Tech. rep., University of
  Illinois at Urbana-Champaign, URL: {\it
  http://theory.cs.uiuc.edu/~jeffe/open/klee.html}.

\bibitem[{Ferrari \& Cribari-Neto(1998)}]{FER98}
Ferrari, S. \& Cribari-Neto, F. (1998) On bootstrap and analytical bias
  correction. \emph{Economics Letters} \textbf{58}, 7--15.

\bibitem[{Figueira \emph{et~al.}(2005)Figueira, Greco \& Erhgott}]{FIG05}
Figueira, J., Greco, S. \& Erhgott, M. (2005) \emph{Multiple criteria decision
  analysis - {S}tate of the art - {S}urvey}. Springer's International Series.

\bibitem[{Fleischer(2003)}]{FLE03}
Fleischer, M. (2003) The measure of {P}areto optima. {A}pplications to
  multi-objective metaheuristics. vol. 262, pp. 519--523, Faro, Portugal.

\bibitem[{Glynn \emph{et~al.}(2009)Glynn, Rubino \& Tuffin}]{GLY09}
Glynn, P., Rubino, G. \& Tuffin, B. (2009) \emph{Robustness properties and
  confidence interval reliability. In: {\it Rare Event Simulation}}. Wiley.

\bibitem[{Hurtado(2004)}]{HUR04}
Hurtado, J. (2004) An examination of methods for approximating implicit limit
  state functions from the viewpoint of statistical learning theory.
  \emph{Structural Safety} \textbf{26}, 271--293.

\bibitem[{Kleijnen(2011)}]{KLE11}
Kleijnen, J. (2011) Simulation optimization via bootstrapped kriging: Survey.
  Tech. rep., Report from Tilburg University, Center for Economic Research.

\bibitem[{Kleijnen \& van Beers(2009)}]{KLE09}
Kleijnen, J. \& van Beers, W. (2009) Monotonicity-preserving bootstrapped
  kriging metamodels for expensive simulations. Tech. rep., Discussion Paper
  2009-75, Tilburg University, Center for Economic Research.

\bibitem[{Kroese \& Rubinstein(2007)}]{KRO07}
Kroese, D. \& Rubinstein, R. (2007) \emph{Simulation and the Monte Carlo Method
  (2nd edition)}. Wiley.

\bibitem[{Lemaire \& Pendola(2006)}]{LEM06}
Lemaire, M. \& Pendola, M. (2006) {PHIMECA-SOFT}. \emph{Structural Safety}
  \textbf{28}, 130--149.

\bibitem[{Limbourg \emph{et~al.}(2010)Limbourg, de~Rocquigny \&
  Andrianov}]{LIM10}
Limbourg, P., de~Rocquigny, E. \& Andrianov, G. (2010) Accelerated uncertainty
  propagation in two-level probabilistic studies under monotony.
  \emph{Reliability Engineering and System Safety} \textbf{95}, 998--1010.

\bibitem[{Lin(1993)}]{LIN93}
Lin, D. (1993) A new class of supersaturated design. \emph{Technometrics}
  \textbf{35}, 28--31.

\bibitem[{MacKay \emph{et~al.}(1979)MacKay, Beckman \& Conover}]{MAC79}
MacKay, M., Beckman, R. \& Conover, W. (1979) A comparison of three methods for
  selecting values of input variables in the analysis of output from a computer
  code. \emph{Technometrics} \textbf{21}, 239--249.

\bibitem[{Madsen \& Ditlevsen(1996)}]{MAD96}
Madsen, H. \& Ditlevsen, O. (1996) \emph{Structural reliability methods}.
  Wiley.

\bibitem[{Meyer(1972)}]{Meyer72}
Meyer, P.A. (1972) \emph{Martingales and stochastic integrals. Lecture Notes in
  Mathematics Vol. 284}. Springer-Verlag.

\bibitem[{Morio(2011)}]{MORIO2011}
Morio, J. (2011) Influence of input pdf parameters of a model on a failure
  probability estimation. \emph{Simulation Modelling Practice and Theory}
  \textbf{19}, 2244--2255.

\bibitem[{Munoz-Muniga \emph{et~al.}(2011)Munoz-Muniga, Garnier, Remy \&
  de~Rocquigny}]{MUN10}
Munoz-Muniga, M., Garnier, J., Remy, E. \& de~Rocquigny, E. (2011) Adaptive
  directional stratification for controlled estimation of the probability of a
  rare event. \emph{Reliability Engineering and System Safety ({\it in press})}
  .

\bibitem[{Overmars \& Yap(1991)}]{OVE91}
Overmars, M.H. \& Yap, C.K. (1991) New upper bounds in {K}lee's measure
  problem. \emph{SIAM Journal of Computing} \textbf{20}, 1034--1045.

\bibitem[{Rajabalinejad \emph{et~al.}(2011)Rajabalinejad, Meester, van Gelder
  \& Vrijling}]{RAJ10}
Rajabalinejad, M., Meester, L., van Gelder, P. \& Vrijling, J. (2011) Dynamic
  bounds coupled with {M}onte {C}arlo simulations. \emph{Reliability
  Engineering and System Safety} \textbf{96}, 278--285.

\bibitem[{Ranjan \emph{et~al.}(2008)Ranjan, Bingham \& Michailidis}]{RAN08}
Ranjan, P., Bingham, D. \& Michailidis, G. (2008) Sequential experiment design
  for contour estimation from complex computer codes. \emph{Technometrics}
  \textbf{50}, 527--541.

\bibitem[{R\"{u}schendorf(2009)}]{RUS09}
R\"{u}schendorf, L. (2009) On the distributional transform, {S}klar's theorem,
  and the empirical copula process. \emph{Journal of Statistical Planning and
  Inference} \textbf{139}, 3921--3927.

\bibitem[{Shamos \& Hoey(1976)}]{SHA76}
Shamos, M.I. \& Hoey, D. (1976) Geometric intersection problems.
  \emph{Proceedings of the 17th IEEE Symposium about the Foundations of
  Computer Science (FOCS '76)} pp. 208--215.

\bibitem[{Tsompanakis \emph{et~al.}(2007)Tsompanakis, Lagaros, Papadrakakis \&
  Frangopol}]{Tsompanakis2007}
Tsompanakis, Y., Lagaros, N., Papadrakakis, M. \& Frangopol, D.e. (2007)
  \emph{Structural design optimization considering uncertainties}. Taylor \&
  Francis.

\bibitem[{van Leeuwen \& Wood(1981)}]{LEW81}
van Leeuwen, J. \& Wood, D. (1981) The measure problem for rectangular ranges
  in $d-$space. \emph{Journal of Algorithms} \textbf{2}, 282--300.

\end{thebibliography}
\bibliographystyle{besjournals}

%
%
%
%

\clearpage
\renewcommand{\appendixname}{Appendix}
\appendix

%
%
%
%

\section{Proofs}\label{proofs}

\begin{proof}[{\bf Proof of Lemma \ref{lemma.0}}]
An infinite uniform sampling on $\U$ provides on the open sets $(\mathring{\U^-},\mathring{\U^+})$  two topologies constituted by the collections of open subsets $\mathring{\U^-_0},\ldots,\mathring{\U^-_n},\ldots$,  and $\mathring{\U^+_0},\ldots,\mathring{\U^+_n},\ldots$. Hence the sequence $({\U}^-_n,{\U}^+_n)$  define two covers (exhaustions) of  $(\U^-,\U^+)$. Then 
\begin{eqnarray*}
U^-_{\infty}=\bigcup_{k=0}^{\infty} U^-_k \ = \ \U^-, \ & \ & \  
U^+_{\infty}=\bigcup_{k=0}^{\infty} U^+_k \ = \ \U^+
\end{eqnarray*}
and 
$\lim_{n\rightarrow\infty} p^-_n   =  \lim_{n\rightarrow\infty} P({\bf X}\in \U^-_n)
 =  P({\bf X}\in U^-_{\infty})  =  \p$ {by inclusion.}
Similarly, $\lim p^+_n  = \p$. Furthermore, given $p^-_0$ and $p^+_0$,  $p^-_n$ and $1-p^+_n$ are ${\cal{F}}_{n-1}-$adapted  submartingales bounded in $\Ll_p$ $\forall p\geq 1$. Then, from Doob's theorem \citep{Meyer72}, the bounds converge almost surely to $\p$. 
\end{proof}


\begin{proof}[{\bf Proof of Proposition \ref{mle.exist}}]
One may write $\ell''_n(\p)  = \sum_{k=1}^n \tilde{\omega}_k\left(\p\right) S_k\left(\p\right)$ with
\begin{eqnarray}
S_k\left(\p\right)  & = & -1 + \left(p_k-\p\right)\tilde{\omega}_k\left(\p\right)\left(2\p - p^-_{k-1} - p^+_{k-1}\right), \label{sk} \\
\nonumber                    & = & -\tilde{\omega}_k(\p)(\p-p_k)^2.
\end{eqnarray}
Hence $\ell''_n(\p)<0$ in $(p^-_{n-1},p^+_{n-1})$. Besides, $\lim_{\p\rightarrow p^-_{n-1}} \ell'_n(\p)=\infty$  and $\lim_{\p\rightarrow p^+_{n-1}} \ell'_n(\p)=-\infty$. Hence, by twice continuity and   differentiability of $\ell_n(\p)$, the mean value theorem implies the existence and unicity of a MLE $\hat{\p}_n$
in $]p^-_{n-1},p^+_{n-1}[$. 
\end{proof}


\begin{proof}[{\bf Proof of Lemma \ref{lemma.expect}}] We shall proceed by induction. 
Since $p^-_0<\p<p^+_0$, (\ref{expect.1}) and (\ref{expect.2}) hold for $n=0$. Denote $\eta_n=1/(\p-p^-_{n})^2$. For $n\geq 1$, it is assumed that $\E[\eta_n]<\infty$. Then
\begin{eqnarray*}
\E\left[\eta_{n+1}\right] & = & \E\left[\eta_n \E\left[1-\s_{\bf x_{n+1}}|{\cal{F}}_{n}\right]\right] + \E\left[\E\left[{\s_{\bf x_{n+1}}}/{\left(\p - p^{-}_n - \Vol^-_{\bf x_{n+1}}\right)^2}|{\cal{F}}_{n}\right]\right]
\end{eqnarray*}
with $\Vol^-_{\bf x_{n+1}}  =  \int_{\U_{n}\cap\U^-} \1_{\{{\bf x}\preceq {\bf x_{n+1}}\}} \ d{\bf x}$ the additive volume of formerly non-dominated failure points in $\U_{n}$ that are now dominated by the failure point ${\bf x_{n+1}}$. By hypothesis, the first term is always finite. Furthermore, with
\begin{eqnarray*}
\Vol^-_{\bf x_{n+1}} & \leq & \sup\limits_{{\bf x_n}\in\U_n\cap\U^-} \int_{\U_{n}\cap\U^-} \1_{\{{\bf x}\preceq {\bf x_n}\}} \ d{\bf x} \ = \ \sup\limits_{{\bf x_n}\in\bar{\Ss}} \int_{\U_{n}\cap\U^-} \1_{\{{\bf x}\preceq {\bf x_n}\}} \ d{\bf x},
\end{eqnarray*}
(\ref{cond001}) implies that $\Vol^-_{\bf x_{n+1}} < \p - \p^-_{n}$. Since $p^-_n=p^-_{n-1}+\Vol^-_{\bf x_{n}}$, etc., one has $\sum_{k=1}^{n+1}\Vol^-_{\bf x_{k}}  <  \p-p^-_0$. Then
\begin{eqnarray*}
\E\left[\E\left[{\s_{\bf x_{n+1}}}/{\left(\p - p^{-}_n - \Vol^-_{\bf x_{n+1}}\right)^2}|{\cal{F}}_{n}\right]\right] & = & \E\left[{\s_{\bf x_{n+1}}}/{\left(\p - p^-_0 -  \sum_{k=1}^{n+1}\Vol^-_{\bf x_{k}} \right)^2}\right] \ < \ \infty.
\end{eqnarray*}
The same rationale applies to $1/(p^+_{n}-\p)^2$, by symmetry, since $p^{+}_{n+1}=p^{+}_n - \Vol^+_{\bf x_{n+1}}$ with $\Vol^-_{\bf x_{n+1}}  =  \int_{\U_{n}\cap\U^+} \1_{\{1-{\bf x}\preceq 1-{\bf x_{n+1}}\}} \ d{\bf x}$.
\end{proof}


\begin{proof}[{\bf Proof of Proposition \ref{asympt.var.reduction}}]
One has $\E[\ell'_n(\p)]  =  \sum_{k=1}^n \E\left[\tilde{\omega}_k(\p)\E\left[p_{k}-\p|{\cal{F}}_{k-1}\right]\right]  =  0$ since $\tilde{\omega}_{n+1}$ depends only on ${\cal{F}}_n$,
hence the Fisher information $J_n(\p)=\Var[\ell'^2_n(\p)]=\E[\ell'^2_n(\p)]$ is equal to $-\E[\ell''_n(\p)]$ by twice differentiability and continuity of $\ell_n(\cdot)$, similarly to a classic iid. case. Assumption 4 implies that $\forall n<\infty$, $p^-_{n-1}<\p<p^{+}_{n-1}$, ie. $\p$ cannot be reached in any finite number of iterations, so that these quantities are well defined. 
With  $-S_n(\p)  =  \tilde{\omega}_n(\p)(\p-p_n)^2$ $\forall~n\geq 0$ from (\ref{sk}),
\begin{eqnarray*}
J_n(\p) 
       & = & \sum\limits_{k=1}^n \E\left[\tilde{\omega}^2_k(\p) \Var\left[p_k|{\cal{F}}_{k-1}\right]\right] \ = \ 
       \sum\limits_{k=1}^n \E\left[\tilde{\omega}_k(\p)\right]
\end{eqnarray*}
since
\begin{eqnarray}
\nonumber \Var\left[p_n|{\cal{F}}_{n-1}\right] & = & \left(p^+_{n-1}-p^-_{n-1}\right)^2\E\left[\s_{\bf x_n}|{\cal{F}}_{n-1}\right] - \left(\p - p^-_{n-1}\right)^2, \\
\nonumber & = & \left(p^+_{n-1}-p^-_{n-1}\right)\left(\p - p^-_{n-1}\right)- \left(\p - p^-_{n-1}\right)^2, \\
& = & \tilde{\omega}^{-1}_n(\p). \label{var.pn}
\end{eqnarray}
Inequality (\ref{exp2.var}) is a simple consequence of Jensen's inequality: since $\E^{-1}\left[\tilde{\omega}^{-1}_k(\p)\right] \leq \E\left[\tilde{\omega}_k(\p)\right]$, then 
${\displaystyle 
J^{-1}_n(\p) \ \leq  \ \left(\sum\limits_{k=1}^n \E^{-1}\left[\tilde{\omega}^{-1}_{k}(\p)\right]\right)^{-1} \ = \   \frac{\p(1-\p)}{\sum\limits_{k=1}^n (1-c_{k-1})^{-1}}.}$
\end{proof}


\begin{proof}[{\bf Proof of Proposition \ref{asympt.var.reduction.2}}]
Using the notation $S_k(\p)$ defined in (\ref{sk}), 
\begin{eqnarray*}
J_n(\p) & = & - {\displaystyle \E\left[\frac{\p(1-\p)}{n V^{MC}_n(\p)}\sum\limits_{k=1}^n \tilde{\omega}_k\left(\p\right) S_k\left(\p\right)\right]} \ = \ n^{-1} \frac{\tilde{J}_n(\p)}{V^{MC}_n(\p)}
\end{eqnarray*}
with 
\begin{eqnarray*} 
\tilde{J}_n(\p) & = & \E\left[\sum\limits_{k=1}^n \p(1-\p)\left(p^+_{k-1}-\p\right)^{2\s_{\bf x_k}-2} \left(\p-p^-_{k-1}\right)^{-2\s_{\bf x_k}}\right],
\end{eqnarray*}
which can be rewritten as
\begin{eqnarray*} 
\tilde{J}_n(\p)    & = & \sum\limits_{k=1}^n \left\{ \E\left[\s_{\bf x_k}\left(\frac{\p(1-\p)}{\left(\p-p^-_{k-1}\right)^2}\right)\right] 
    + 
    \E\left[(1-\s_{\bf x_k})\left(\frac{\p(1-\p)}{\left(p^+_{k-1}-\p\right)^2}\right)\right]\right\}.
\end{eqnarray*}
Since
 ${\p}^{-1}\left(\left(p^+_{k-1}-\p\right) + \left(\p-p^-_{k-1}\right)\right)  \leq  \rho_{k-1}$, then
\begin{eqnarray*}
p^+_{k-1}-\p & \leq & \p(\rho_{k-1}-1)+p^-_{k-1} \ \leq \ \p(\rho_{k-1}-1)+\p = \p\rho_{k-1}, \\  
\p-p^-_{k-1} & \leq & \p(\rho_{k-1}+1)-p^+_{k-1} \ \leq \ \p(\rho_{k-1}+1) -\p = \p\rho_{k-1}.
\end{eqnarray*}
Hence
\begin{eqnarray*}
\frac{\p(1-\p)}{\left(p^+_{k-1}-\p\right)^2}  \geq  \frac{1-\p}{\p\rho^2_{k-1}} & \text{and} & 
\frac{\p(1-\p)}{\left(\p-p^-_{k-1}\right)^2}  \geq \frac{1-\p}{\p\rho^2_{k-1}}. \\
\end{eqnarray*}
Consequently,
\begin{eqnarray*}
\tilde{J}_n(\p) & \geq & {\displaystyle \frac{1-\p}{\p} \sum\limits_{k=1}^n \E\left[\frac{1}{\rho^2_{k-1}}\left(\s_{\bf x_k} +  1-\s_{\bf x_k}\right)\right]}, \\
                & \geq & {\displaystyle \frac{1-\p}{\p} n \E\left[\frac{1}{\rho^2_{0}}\right]}
\end{eqnarray*}
since $(\rho_n)$ is a strictly decreasing positive process. Since $(p^-_n,p^+_n)$ are predictible processes, $p^-_0$ and $p^+_0$ are deterministic quantities, then 
\begin{eqnarray*}
\E\left[\frac{1}{\rho^2_{0}}\right] & = & \left(\frac{p^-_0}{p^+_0-p^-_0}\right)^2 \ = \ 1/\gamma_0,
\end{eqnarray*}
and $\tilde{J}_n(\p)\geq n\left(\frac{1-\p}{\gamma_0\p}\right)$ which proves (\ref{exp1.var}). 
\end{proof}


\begin{proof}[{\bf Proof of Theorem \ref{asympt.norm.MLE}}]
Given the strong consistency of $\hat{\p}_n$, its asymptotic normality can be established using arguments studied by \cite{CRO75,CRO83}. Showing that $\ell'_n(\p)$ is a ${\cal{F}}_{n-1}-$adapted martingale is a classic result: 
\begin{eqnarray*}
\E\left[\ell'_{n+1}(\p)-\ell'_{n}(\p)|{\cal{F}}_n\right] & = & \tilde{\omega}_{n+1}(\p)\E\left[p_{n+1}-\p|{\cal{F}}_n\right] \ = \ 0.
\end{eqnarray*}
Furthermore $J_n(\p)<n \E[\tilde{\omega}_n(\p)]<\infty$ under Assumption 4 (cf. Lemma \ref{lemma.expect}) Hence $\ell'_n(\p)$ is square integrable. Denoting $\Delta_n(\p)=\ell'_n(\p) - \ell'_{n-1}(\p)$, then $\Delta^2_n(\p)  =  \tilde{\omega}^2_n(\p)(p_n-\p)^2$ and
\begin{eqnarray*}
\E\left[\Delta^2_n(\p)|{\cal{F}}_{n-1}\right] & = & \tilde{\omega}^2_n(\p) \Var\left[p_n|{\cal{F}}_{n-1}\right] \ = \ \tilde{\omega}_n(\p).
\end{eqnarray*}
Then $<\ell'(p)>_{{}_n}  = \sum_{k=1}^n \tilde{\omega}_k(\p)$ denotes the increasing (or bracket) process of $\ell'_n(\p)$. The proof can be achieved in three steps. 

\paragraph*{\bf 1} With $J_n(\p)=\E[<\ell'(p)>_{{}_n}]$ and $\lim_{n\to\infty} J_n(\p)=\infty$ from (\ref{exp1.var}), establishing asymptotic normality first requires to prove the following law of large numbers (LLN)
\begin{eqnarray}
M_n(\p) \ = \ J^{-1}_n(\p)<\ell'(p)>_{{}_n} - 1 & \xrightarrow[]{\Pp} & 0. \label{MLLN}
\end{eqnarray}
Denote  $W_n(\p)=\sum_{k=1}^n (\tilde{\omega}_{k}(\p)-\E[\tilde{\omega}_{k}(\p)])$ Then, $\forall \epsilon > 0$,
\begin{eqnarray*}
P\left(|M_n(\p)|>\epsilon\right) & = & P\left(J^{-1}_n(\p)|W_n(\p)|>\epsilon\right), \\
 & \leq & P\left(V^{MC}_n(\p) \frac{\p\gamma_0}{1-\p}|W_n(\p)|>\epsilon\right)\ \ \ \text{from (\ref{exp1.var})}, \\
 & \leq & P\left(\frac{1}{n}|W_n(\p)|>\epsilon'\right)\ \ \ \text{ with $\epsilon'=\epsilon/(\p^2\gamma_0)$,}
\end{eqnarray*}
which tends to 0 under {\bf (i)} and proves (\ref{MLLN}). \\

\paragraph*{\bf 2} For all $k\in\{1,\ldots,n\}$, denote $\Gamma_{k,n} = J^{-1/2}_n(\p){|\Delta_k(\p)|}$. The second requirement of asymptotic normality is proving the following Lindeberg condition: $\forall \epsilon>0$, 
\begin{eqnarray}
\frac{1}{J_n(\p)}\sum\limits_{k=1}^n \E\left[\Delta^2_k(\p)\1_{\left\{\Gamma_{k,n}> \epsilon\right\}}|{\cal{F}}_{k-1}\right] & \xrightarrow[n\to\infty]{\Pp} & 0. \label{lindeberg.10}
\end{eqnarray} 
A Lyapunov condition is often used instead of (\ref{lindeberg.10}), but  requires $2+\delta$-order moment assumptions on $\tilde{\omega}_k(\p)$. An alternative approach is the following. From Markov's inequality and since the $\tilde{\omega}_k(\p)$ are increasing functions of $k$,
\begin{eqnarray*}
P\left(\Gamma_{k,n}>\epsilon|{\cal{F}}_{k-1}\right) & \leq & \frac{\tilde{\omega}_k(\p)}{\epsilon^2 J_n(\p)} 
 \  \leq \  \frac{\tilde{\omega}_n(\p)}{\epsilon^2 J_n(\p)}. 
\end{eqnarray*}
It follows from (\ref{MLLN}) that 
\begin{eqnarray*}
\frac{\tilde{\omega}_n(\p)}{J_n(\p)} + \frac{1}{J_n(\p)}\sum\limits_{k=1}^{n-1}\tilde{\omega}_n(\p) & = & \frac{\tilde{\omega}_n(\p)}{J_n(\p)} + \left(\frac{J_{n-1}(\p)}{J_n(\p)}\right)\left(\frac{<\ell'(p)>_{{}_{n-1}}}{J_{n-1}(\p)}\right), \\
& \xrightarrow[n\to\infty]{\Pp} & 1.
\end{eqnarray*} 
However, by Lemma \ref{lemma.expect}, $\E[\tilde{\omega}_n(\p)]<\infty$ which means that $J_n(\p)\overset{\infty}{\sim} J_{n-1}(\p)$. Necessarily, ${\tilde{\omega}_n(\p)}/{J_n(\p)}  \xrightarrow[n\to\infty]{\Pp}  0$ 
and 
\begin{eqnarray}
L_{k,n} \ = \ \E\left[\1_{\left\{\Gamma_{k,n}> \epsilon\right\}}|{\cal{F}}_{k-1}\right] & \xrightarrow[k \to n \to\infty]{\Pp} & 0. \label{LNN30}
\end{eqnarray}
Note besides that
\begin{eqnarray*}
\E\left[\Delta^2_k(\p)\1_{\left\{\Gamma_{k,n}> \epsilon\right\}}|{\cal{F}}_{k-1}\right] & \leq & \E\left[\Delta^2_k(\p)|{\cal{F}}_{k-1}\right]\E\left[\1_{\left\{\Gamma_{k,n}> \epsilon\right\}}|{\cal{F}}_{k-1}\right] + \left|\mbox{Cov}\left[\Delta^2_k(\p),\1_{\left\{\Gamma_{k,n}> \epsilon\right\}}|{\cal{F}}_{k-1}\right]\right|, \\
& \leq & \tilde{\omega}_k(\p) L_{k,n} + \tilde{\omega}^2_k(\p) \sqrt{\Var\left[(p_k-p)^2|{\cal{F}}_{k-1}\right]}\sqrt{\Var\left[\1_{\left\{\Gamma_{k,n}> \epsilon\right\}}|{\cal{F}}_{k-1}\right]}
\end{eqnarray*}
from Cauchy-Schwarz inequality. Since $\Var[X^2]\leq \E[X]$ when $X\in\{0,1\}$, then
$\sqrt{\Var[\1_{\{\Gamma_{k,n}> \epsilon\}}|{\cal{F}}_{k-1}]}  \leq  \sqrt{L_{k,n}}$.
Furthermore, denote
\begin{eqnarray*}
K_{k,n}(\p) & = & \tilde{\omega}_k(\p) \sqrt{\Var\left[(p_k-p)^2|{\cal{F}}_{k-1}\right]}. 
\end{eqnarray*}
From Lemma \ref{lemma.id} below and under {\bf (ii)}, then $K_{k,n}(\p) \xrightarrow{\Pp}{} 0$. 
Therefore, one may write
\begin{eqnarray*}
\E\left[\Delta^2_k(\p)\1_{\left\{\Gamma_{k,n}> \epsilon\right\}}|{\cal{F}}_{k-1}\right] & \leq & \tilde{\omega}_k(\p) \beta_{k,n}
\end{eqnarray*}
with $\beta_{k,n} = L_{k,n} + K_{k,n}(\p)\sqrt{L_{k,n}} \xrightarrow{\Pp}{} 0$ from (\ref{LNN30}). Then
\begin{eqnarray}
\frac{1}{J_n(\p)}\sum\limits_{k=1}^n \E\left[\Delta^2_k(\p)\1_{\left\{\Gamma_{k,n}> \epsilon\right\}}|{\cal{F}}_{k-1}\right] & \leq & 
\frac{<\ell'(\p)>_n}{J_n(\p)} \frac{\sum\limits_{k=1}^n \tilde{\omega}_k(\p) \beta_{k,n}}{\sum\limits_{k=1}^n \tilde{\omega}_k(\p)}
\end{eqnarray}
and given (\ref{MLLN}), Toeplitz lemma proves (\ref{lindeberg.10}). Finally,  (\ref{MLLN}) and (\ref{lindeberg.10}) prove the two martingale central limit theorems \citep{BER08}:
\begin{eqnarray}
J^{-1/2}_n(p) \ell'_n(\p) & \xrightarrow[n\to\infty]{{\cal{L}}} & {\cal{N}}(0,1), \label{tlc.martingale.2} \\ 
\frac{\sqrt{J_n(\p)}}{<\ell'(\p)>_{{}_n}} \ell'_n(\p) & \xrightarrow[n\to\infty]{{\cal{L}}} & {\cal{N}}(0,1). \label{tlc.martingale.1} 
\end{eqnarray} 

{
\hbox{\raisebox{0.4em}{\vrule depth 0pt height 0.4pt width 13.25cm} }
\begin{lemma}\label{lemma.id} 
If $\exists \ \gamma_{\infty}$ such that $0<\gamma_{\infty}<\infty$ and $\frac{p^+_n-\p}{\p - p^-_n} \xrightarrow{\Pp}{} \gamma_{\infty}$, then, $\forall k\geq 1$, 
\begin{eqnarray*}
\Var\left[\tilde{\omega}_k(p) (p_k-\p)^2|{\cal{F}}_{k-1}\right] & \xrightarrow{\Pp}{} & \gamma_{\infty} + 1/\gamma_{\infty} -2. 
\end{eqnarray*}
\end{lemma}
\paragraph*{\bf Proof} {\it One may write 
\begin{eqnarray*}
\tilde{\omega}_k(p) (p_k-\p)^2 & = & (1-\s_{\bf x_k})\frac{\p-p^-_{k-1}}{p^+_{k-1}-\p} + \s_{\bf x_k}\frac{p^+_{k-1}-\p}{\p-p^-_{k-1}}, \\
                             & = & \s_{\bf x_k}\left[\tilde{\omega}_k(p)\left\{ \left(p^+_{k-1}-\p\right)^2 - \left(\p-p^-_{k-1}\right)^2\right\}\right] + \frac{\p-p^{-}_{k-1}}{p^+_{k-1}-\p} 
\end{eqnarray*}
With $\s_{\bf x_k}\sim {\cal{B}}(\gamma_k)$ and from (\ref{bernoulli}), then 
\begin{eqnarray*}
\Var\left[\tilde{\omega}_k(p) (p_k-\p)^2|{\cal{F}}_{k-1}\right] & = & \frac{\tilde{\omega}_k(p)}{(p^{+}_{k-1} - p^{-}_{k-1})^2}\left[(p^+_{k-1} + p^-_{k-1} - 2\p)(p^{+}_{k-1} - p^{-}_{k-1})\right]^2, \\
                                                                & = & \tilde{\omega}_k(p)(p^+_{k-1} + p^-_{k-1} - 2\p)^2, \\
                                                                & = & \frac{p^+_{k-1}-\p}{\p-p^-_{k-1}} + \frac{\p-p^-_{k-1}}{p^+_{k-1}-\p} - 2 \ \xrightarrow{\Pp}{} \ \gamma_{\infty} + 1/\gamma_{\infty} -2.
\end{eqnarray*}
\hbox{\raisebox{0.4em}{\vrule depth 0pt height 0.4pt width 13.25cm} }
} 
}

\vspace{1cm}

\paragraph*{\bf 3} A last condition is required to transfer the asymptotic normality from $\ell'_n(\p)$ to $(\hat{\p}_n-\p)$. Since $p^-_{n-1}<\hat{\p}_n<p^{+}_{n-1}$, for any $n$ there always exists an open neighborhood ${\cal{V}}_{\hat{\p}_n}$ of $\p$ containing $\hat{\p}_n$. From twice differentiability of $\ell_n(\cdot)$ and continuity of $\ell'_n(\cdot)$ in ${\cal{V}}_{\hat{\p}_n}$, the mean value theorem  implies there exists some intermediate point $\bar{p}_n\in{\cal{V}}_{\hat{\p}_n}$ between $\p$ and $\hat{\p}_n$ such that
\begin{eqnarray*}
\ell'_n\left(\hat{\p}_n\right) \ = \ 0 \ = \ \ell'_n(\p) + \left(\hat{\p}_n-\p\right)\ell''_n(\bar{p}_n)
\end{eqnarray*}
and moreover $\bar{p}_n\xrightarrow[]{a.s.} \p$. 
Thus, with $\ell''_n(\tilde{p}_n)\neq 0$,
\begin{eqnarray}
\left(\hat{\p}_n-\p\right) & = & \ell'_n(\p)\left(-\ell''_n(\bar{p}_n)\right)^{-1} \label{equ.prov.0}
\end{eqnarray}
and it is necessary to prove the LLN
\begin{eqnarray}
 \frac{\ell''_n(\bar{p}_n)}{<\ell'(\p)>_n} & \xrightarrow[]{\Pp} & 1 \label{LLN40}
\end{eqnarray}
to obtain the final result (Theorem 3 in \citet{CRO83}), combining (\ref{LLN40}) with (\ref{MLLN}) and (\ref{tlc.martingale.2}). Based on {\bf (iii)} this last LLN is straightforward. Indeed, $\forall k\leq n$, 
\begin{eqnarray*}
\left|\frac{p^+_k -\bar{p}_n}{p^+_k - \p} - 1\right| & = & \frac{|\bar{p}_n - \p|}{p^+_k - \p}
\ \leq  \   \frac{|\bar{p}_n - \p|}{p^+_n - \p} \ \xrightarrow[]{\Pp} \ 0.
\end{eqnarray*}
Similarly $({\bar{p}_n-p^-_k})/({\p-p^-_k}) - 1 \xrightarrow[k \to n \to \infty]{\Pp} \ 0$. With $p_{k+1}\in\{p^-_k,p^+_k\}$ then, for $k\in\{0,\ldots,n-1\}$, 
\begin{eqnarray*}
\gamma_{k+1,n} & = & \left(\bar{p}_n - p_{k+1}\right)^2\tilde{\omega}_{k+1}(\bar{p}_n) \ \xrightarrow[k \to n \to \infty]{\Pp} \ 1.
\end{eqnarray*}
Furthermore, some calculus proves that $\kappa_{k,n} = \tilde{\omega}_{k+1}(\bar{p}_n)/\tilde{\omega}_{k+1}(\p)  \xrightarrow[k \to n \to \infty]{\Pp}  1$. Then,
with $-\ell''_n(\p)=\sum_{k=1}^n \Delta^2_k(\p)$, 
\begin{eqnarray*}
 \frac{\ell''_n(\bar{p}_n)}{<\ell'(\p)>_n} & = & \frac{\sum\limits_{k=1}^n \tilde{\omega}_k(\p) \kappa_{k,n} \gamma_{k,n}}{\sum\limits_{k=1}^n \tilde{\omega}_k(\p)} \ \xrightarrow[k \to n \to \infty]{\Pp} \ 1 \ \ \ \text{from Toeplitz lemma.}
\end{eqnarray*}
\end{proof}


\begin{proof}[{\bf Proof of Proposition \ref{TCL.variance}}]
Using the notations of the previous proof, note that $\hat{J}_n(\p)=<\ell'(\p)>_n$. By twice continuity and derivability of $\hat{J}^{~-1}_n(.)$ in $]p^-_n,p^+_n[$, a Taylor expansion gives 
\begin{eqnarray*}
\hat{J}^{~-1}_n(\hat{\p}_n) & = & \hat{J}^{~-1}_n(\p) - \frac{\hat{J}'_n(\p)}{\hat{J}_n^{~2}(\p)}\left(\hat{\p}_n-\p\right)\left(1 + o(1)\right).
\end{eqnarray*}
After some  calculus,
\begin{eqnarray*}
\frac{\hat{J}_n^{~5/2}(\p)}{|\hat{J}'_n(\p)|}\left(\hat{J}^{~-1}_n(\hat{\p}_n)-J^{-1}_n(\p)\right) & = &  R_n U_n + R_n Z_n
\end{eqnarray*}
with $R_n = \sqrt{J_n(\p)/\hat{J}_n(\p)}  \xrightarrow[]{\Pp} 1$ from (\ref{MLLN}), $U_n  =  \mbox{sgn}(\hat{J}'_n(\p))\sqrt{J_n(\p)}\left(\hat{\p}_n-\p\right)(1+o(1))\xrightarrow[]{{\cal{L}}} {\cal{N}}(0,1)$ from Theorem \ref{asympt.norm.MLE}, and 
\begin{eqnarray*}
Z_n & = & \frac{\hat{J}_n^{~3/2}(\p)}{|\hat{J}'_n(\p)|}\left(\frac{\hat{J}_n(\p)}{J_n(\p)} - 1\right). 
\end{eqnarray*}
Thanks to Slutsky's theorem, it is enough to show that $Z_n\xrightarrow[]{\Pp} 0$ to
prove the statement of the proposition. Notice that
\begin{eqnarray*}
\hat{J}'_n(\p) & = & \sum\limits_{k=1}^n \tilde{\omega}^2_k(\p) \left\{2\p - \left(p^+_{k-1} + p^-_{k-1}\right)\right\}
\end{eqnarray*}
which is always nonzero assuming {\bf (iv)}. 
 H{\"o}lder's inequality gives
\begin{eqnarray*}
\frac{\sum\limits_{k=1}^n \tilde{\omega}_k(\p)}{\sum\limits \tilde{\omega}^2_k(\p) \left\{2\p - \left(p^+_{k-1} + p^-_{k-1}\right)\right\}} & \leq & \frac{\sum\limits \left\{2\p - \left(p^+_{k-1} + p^-_{k-1}\right)\right\}^{-1}}{\sum\limits_{k=1}^n \tilde{\omega}_k(\p)}
\end{eqnarray*}
hence
\begin{eqnarray*}
\frac{\hat{J}_n^{~3/2}(\p)}{|\hat{J}'_n(\p)|} & \leq & \frac{\sum\limits \left\{2\p - \left(p^+_{k-1} + p^-_{k-1}\right)\right\}^{-1}}{\sqrt{\sum\limits_{k=1}^n \tilde{\omega}_k(\p)}}
\end{eqnarray*}
Another H{\"o}lder's inequality gives
${\displaystyle \frac{\hat{J}_n^{~3/2}(\p)}{|\hat{J}'_n(\p)|}  \leq  \sqrt{\sum_{k=1}^n \frac{\left(p^+_{k-1}-\p\right)\left(\p-p^-_{k-1}\right)}{2\p - \left(p^+_{k-1} + p^-_{k-1}\right)}}}$
and simple calculus shows that each term of the sum  is stricly smaller than 1. Then
\begin{eqnarray*}
|Z_n| & \leq & \sqrt{n}\left|\frac{\hat{J}_n(\p)}{J_n(\p)} - 1\right| \ \leq \ \frac{2\p^2}{\sqrt{n}} \sum\limits_{k=1}^n \left(\tilde{\omega}_k(\p) - \E\left[\tilde{\omega}_k(\p)\right] \right) 
\end{eqnarray*}
 from (\ref{exp1.var}), then $Z_n\xrightarrow[]{\Pp} 0$ if {\bf (i)} remains true $\forall \delta\geq 1/2$.  
 \end{proof}

\clearpage

%
%
%
%

\section{Supplementary Material}

This supplementary section first provides details about the implementation of sweepline algorithms to solve Klee's measure problem, which allows for an exact computation (modulo rounding errors) of the probability bounds $(p^-_n,p^+_n)$ ; a pseudo-code is given for direct use. Then 
a general result is given about the preservation of monotonicity when the uniform input ${\bf x}=(x_1,\ldots,x_d)$ results from an inverse transformation of the joint cdf.  

\vspace{-1.5cm}

\Subsection{A sweepline  algorithm to compute volumes of hypercubic unions}

Sweepline (or {\it plane sweep}) algorithms are commonly used to
jointly detect and sort intersections between segments \citep{LEW81}. The
$d$-dimensional volume is calculated recursively by exploring 
all n-1-dimensional ``slices" of the $d$-th dimension. See
\citet{SHA76,deB97} and \citet{CHL98} for more explanations. When segments are
parallel or perpendicular such as their intersections define a
union of hypercubes sharing the same orthogonal basis, the volume calculation is known as Klee's
measure problem \citep{ERI98,CHA08}.  A pseudo-code follows to be used for
direct implementation. 

Let $\Delta_n$ be the $n\times d$ matrix of
$n$ vertexes $({\bf x}_1,\ldots,{\bf x}_n)$ defining the union of
hypercubes (for an example, see Figure \ref{graphic.example.2}).
In the following pseudo-code, the volume considered is $V^-_n$,
also defined
by the points of $\Delta_n$ and the origin $(0,\ldots,0)$ of the $\U-$space. \\

\texttt{ \hspace{-0.5cm}\noindent{\bf Algorithm}
$\mbox{VOL}(\Delta_n,n,d){\bf .}$ 
\hbox{\raisebox{0.4em}{\vrule depth 0pt height 0.4pt width 10cm} }
\begin{enumerate}
\item Let $\Delta'_n = \sigma_{n,d}(\Delta_n)$ be the $n\times d$ permutation of $\Delta_n$ arranged in the increasing order
of the $n-$vector of $d-$dimensional components.
\item Remove the $d-$dimensional components from $\Delta'_n$ and denote $\Vol_n=0$.
\item For $i\in\{1,\ldots,n\}$,
\begin{enumerate}
\item Consider the slice $\Delta^{(i)}_n = \left\{{\bf x}'_{i},\ldots,{\bf x}'_n \in \Delta'_n\right\}$.
\item Denote $\widetilde{\Vol}^{(i)}_{n}$ the $d-1-$dimensional volume of $\Delta^{(i)}_n$. 
\begin{itemize}
\item[{}] If $\dim Z^{(i)}_n = 1$, 
\begin{itemize} 
 \item[\textbullet] $\Delta^{(i)}_n$ is a $n-i+1-$vector and
$\widetilde{\Vol}^{(i)}_{n} = \max\{{\bf x}\in \Delta^{(i)}_n\}$;
 \item[\textbullet]  force $i$ to the index of this maximal component in $\Delta'_n$; 
\end{itemize}
\item[{}] else $\widetilde{\Vol}^{(i)}_{n} = \mbox{VOL}(\Delta^{(i)}_n,n-i+1,d-1)$.
\end{itemize}
\vspace{0.15cm} \item Let $\Lambda_i = \Delta'_n[i,d] -
\Delta'_n[i-1,d]$ the size of $\Delta^{(i)}_n$ (assuming
$\Delta'_n[0,d]=0$). \item Compute $\Vol^{(i)}_{n} = \Lambda_i
\cdot\widetilde{\Vol}^{(i)}_{n}$ the $d-$dimensional  volume.
\item Update the total volume $\Vol_n = \Vol_n + \Vol^{(i)}_{n}$.
\end{enumerate}
\end{enumerate}
\hbox{\raisebox{0.4em}{\vrule depth 0pt height 0.4pt width 14.5cm} }
}

In practice, this algorithm seems to remain little used for dimension $d$ larger than 2 or 3. This is not surprising  because its complexity $C(n,d)$ (the number of runs for a $d$-dimensional hypervolume between $n$ points) is $O(n^d)$. This appears when considering the first developments of $C(n,d)$:
\begin{eqnarray*}
C(n,d) & = & \sum\limits_{k=0}^{n-1} C(n-k,d-1) \ = \  \sum\limits_{k=0}^{n-1} (k+1) C(n-k,d-2), \\
       & = & \sum\limits_{k=0}^{n-1} \left(\sum\limits_{p=0}^{k-1} p\right) C(n-k,d-3), \\
       &  = & \sum\limits_{k=0}^{n-1} \frac{(k+1)(k+2)}{2} \ C(n-k,d-3), \\
       & \ldots &
\end{eqnarray*}
Note however than the fastest version of this algorithm,  proposed by 
\citet{OVE91}, runs in time $O(n^{d/2}\log n)$ for $d\geq 3$. An alternative approach was presented 
by \citet{CHL98} with the same asymptotic performance, although its exposition was restricted to
dimensions 3 and 4. At the present time the computational difficulties raised by diminishing the cost still remain open problems, although some slight improvements have recently been found by \citet{CHA08}. Some ideas of possible improvements could possibly come from a parallel with multi-objective optimization contexts (cf. Remark 1 in the article). Indeed, algorithms running in polynomial time $O(n^{k_1} d^{k_2})$ to compute hypervolume metrics of Pareto frontiers have already been proposed by \citet{FLE03}. 

\Subsection{Preservation of monotonicity through space transformation}\label{space.transformation.appendix}

Consider $\tilde{g}$ a monotonic function with physical input random vector ${\bf y}=(y_1,\ldots,y_d)$ and denote $T$ their multivariate distributional transform. The methodology proposed in the article applies using the transformed function $g=\tilde{g}\circ T^{-1}$, provided $T^{-1}$ is a globally increasing function of independent uniform inputs ${\bf x}=(x_1,\ldots,x_d)$. This is ensured when $(y_1,\ldots,y_d)$ are independent, since $T^{-1}=(F^{-1}_1,\ldots,F^{-1}_d)$ where $F_i$ is the $i$th marginal cumulative distribution function (cdf). In dependent cases (and possibly when the physical inputs mix continuous and discrete distributions), the generalized Rosenblatt's transform \citep{RUS09} may be used if the inputs can be stochastically conditioned, namely they can be sorted to get the explicit writing 
\begin{eqnarray*}
T(y_1,\ldots,y_d)  & =  & F_1(y_1)\prod_{i=2}^d F_{i|1,\ldots,i-1}(y_i|y_1,\ldots,y_{i-1}).
\end{eqnarray*}
Under this assumption, next lemma provides an intuitive sufficient condition for $F^{-1}$ to be an increasing function of all $x_i\sim{\cal{U}}[0,1]$. 

\begin{lemma}\label{sufficient.cond}
Assume that for $i=2,\ldots,d$, there exists a mapping $f_i$ and a set of (possibly random) parameters $\theta_i$ independent of $Y_1,\ldots,Y_i$ such that:
\begin{description}
\item[(i)] $Y_i=f_i(Y_1,\ldots,Y_{i-1},\theta_i)$,
\item[(ii)] $f_i$ is a globally increasing function of $Y_1,\ldots,Y_{i-1}$;
\end{description}
then 
 $T^{-1}({\bf x})$ is an increasing function of ${\bf x}$. 
\end{lemma}

Multivariate normal distributions are often selected as approximate ways to tackle the difficulties of assessing correlations between input physical parameters, and therefore deserve particular interest in the field of computer experiments. If \citet{CHEN09} obtained general results about the preservation of monotonicity when these distributions are given under the form of Gaussian copulas, an immediate corollary of Lemma \ref{sufficient.cond} is to notice that any standard binormal input distribution with positive correlation coefficient $\mu$  ensures that $T^{-1}({{\bf x}})$ is increasing. Indeed, ${\bf Y}=(Y_1,Y_2)$ where $Y_1  \sim  {\cal{N}}(0,1)$ and $Y_2  =  \mu Y_1 + \sqrt{1-\mu^2} \theta$ with $\theta\sim{\cal{N}}(0,1)$. A similar result can be found for the class of elliptical  bivariate copulas. \\

\begin{proof}[{\bf Proof of Lemma \ref{sufficient.cond}}]
Assume {\bf (i)}. $\forall t\in\R, \ \ \forall k\in\{2,\ldots,d\}$, denote
$ p^t_{\theta_i}\left(Y_1,\ldots,Y_{i-1}\right)  =  P\left(f_i(Y_1,\ldots,Y_{i-1},\theta_i)< t|Y_1,\ldots,Y_{i-1}\right)$.
Then, $\forall z\in\R$, let $A^t_{Y_1,\ldots,Y_{i-1}}(z)$ denote the event $\{p^t_{\theta_i}\left(Y_1,\ldots,Y_{i-1}\right)\leq z\}$. 
By definition,
\begin{eqnarray*}
F^{-1}_{i|1,\ldots,i-1}\left(z|Y_1,\ldots,Y_{i-1}\right) & = & \inf\left\{t \in \R \ | {{P}}\left(A^t_{Y_1,\ldots,Y_{i-1}}(z)\right)=1\right\}.
\end{eqnarray*}
Assuming {$(ii)$}, $p^t_{\theta_i}\left(Y_1,\ldots,Y_{i-1}\right)$ is a globally decreasing function of $Y_1,\ldots,Y_{i-1}$. Thus, given $t$, the occurence of event $A^t_{Y_1,\ldots,Y_{i-1}}(y)$ similarly decreases. Necessarily $t$ increases, hence the minimum value of all $t\in\R$ such that ${\cal{P}}(A^t_{Y_1,\ldots,Y_{i-1}}(z))=1$ increases. Hence $F^{-1}_{i|1,\ldots,i-1}$ is a globally increasing function of $Y_1,\ldots,Y_{i-1}$, $\forall i\in\{2,\ldots,d\}$. Since $Y_1=F^{-1}(X_1)$ is naturally an increasing function of $X_1$, a simple recursive reasoning shows that $F^{-1}_{i|1,\ldots,i-1}$ is an increasing function of $X_1,\ldots,X_{i-1}$. The statement of the lemma follows.
\end{proof}

\end{document}